\documentclass[10pt]{amsart}
\usepackage{amsmath}
\usepackage{amsthm}
\usepackage{amsfonts}
\usepackage{amssymb}
\usepackage{latexsym}
\usepackage{enumerate}

\def\const{{\rm const}}

\newcommand{\qbar}{\underline{q}}
\newcommand{\beq}{\begin{eqnarray}}
\newcommand{\eeq}{\end{eqnarray}}
\newcommand{\beqno}{\begin{eqnarray*}}
\newcommand{\eeqno}{\end{eqnarray*}}
\newcommand{\bR}{\mathbb R}
\newcommand{\bN}{\mathbb N}
\newcommand{\cN}{\mathcal N}
\newcommand{\cM}{\mathcal M}
\newcommand{\cG}{\mathcal G}

\newcommand{\cF}{\mathcal F}

\newcommand{\cX}{\mathcal X}

\newcommand {\Bh}{B^{(h)}}
\newcommand {\Bvh}{B^{(h)}_{4}}

\newcommand{\midint}{- \mskip-18mu \int}

\newcommand{\eps}{\varepsilon}
\renewcommand{\phi}{\varphi}

\DeclareMathOperator{\spt}{supp}

\allowdisplaybreaks

\def\mvint_#1{\mathchoice
  {\mathop{\vrule width 6pt height 3 pt depth -2.5pt
        \kern -8pt \intop}\nolimits_{\kern -3pt #1}}%
 {\mathop{\vrule width 5pt height 3 pt depth -2.6pt
                  \kern -6pt \intop}\nolimits_{#1}}%
 {\mathop{\vrule width 5pt height 3 pt depth -2.6pt
                  \kern -6pt \intop}\nolimits_{#1}}%
 {\mathop{\vrule width 5pt height 3 pt depth -2.6pt
                  \kern -6pt \intop}\nolimits_{#1}}}

\newtheorem{theorem}{Theorem}[section]

\newtheorem{Lem}[theorem]{Lemma}

\newtheorem*{Bem}{Remark}

\theoremstyle{definition}

\newtheorem*{remark}{Remark}

\numberwithin{equation}{section}

\begin{document}

\title{Partial regularity for minima of higher order functionals with $p(x)$ growth}

\author{Jens Habermann}

\date{\scriptsize \today}

\begin{abstract}
For higher order functionals $\int_\Omega f(x,\delta u(x), D^m u(x))\, dx$ with 
$p(x)$-- growth with respect to the variable containing $D^m u$, we prove 
that $D^m u$ is H\"older continuous on an open subset $\Omega_0 \subset \Omega$ of
full Lebesgue-- measure, provided that the exponent function $p:\Omega \to (1,\infty)$ itself is H\"older continuous.
\end{abstract}

\maketitle

\section{Introduction}

In this paper we are concerned with regularity for minimizers of quasiconvex 
functionals of higher order with $p(x)$-- growth. 

We consider functionals of the type
\begin{equation}\label{varint}
     \cF\left[w,\Omega\right] = \int_{\Omega} f\left(x,\delta w(x), 
     D^m w(x)\right) \, dx, 
\end{equation}
on the space $W^{m,1}_{loc}(\Omega;\bR^N),\ N > 1$, where $\Omega \subset \bR^n$ is an open bounded domain and $f: \Omega \times 
\bR^N \times \bR^{nN} \times \ldots \times \bR^{N \binom{n+m-1}{m}} \to \bR$ a 
Carath\'eodory function. $\delta w \equiv (w,Dw,\ldots,D^{m-1}w)$ denotes the
vector containing the lower order derivatives. For $k = 1,\ldots,m$ we use the
notation $D^k u \equiv \left\{ D^\alpha u_i\right\}_{i=1,\ldots,N}^{|\alpha| = k}$
for the derivative of order $k$. Note that $D^k u$ is an element of the space 
$\odot^k (\bR^n;\bR^N)$ of
symmetric $k$-- linear forms on $\bR^n$ with values in $\bR^N$ which can be identified with the space $\bR^{N \binom{n+k-1}{k}}$. In the whole paper, for the seek of brevity we use the
abbreviations ${\mathcal M} \equiv N\sum_{k=0}^{m-1} \binom{n+k-1}{k}$ and
${\mathcal N} \equiv N\binom{n+m-1}{m}$. With this notation we have $D^m u(x) \in
\bR^{\mathcal N}$ and $\delta u(x) \in \bR^{\mathcal M}$.

The functional above is supposed to have $p(x)$-- growth, i.e. $p: \Omega \to
(1,\infty)$ is a continuous function and 
\begin{equation}\label{Intr.Wachst.}
f(x,\xi,z) \approx \left( 1 + |z|^2 \right)^{p(x)/2},
\end{equation}
for all $x \in \Omega$, $\xi \in \bR^{\cM}$ and $z \in \bR^{\cN}$.
Additionally we suppose the functional to be uniformly strictly quasiconvex (see Section \ref{sec:setting} for the exact definition).

We call $u \in W_{loc}^{m,1} \left( \Omega , \bR^N \right)$ a {\em local minimizer} of
the functional $\cF$, if 
\begin{equation}\label{p(x)-Integr.}
     \left| D^m u \right|^{p(\cdot)} \in L^1_{loc} \left( \Omega \right),
\end{equation}
and
\begin{equation}\label{Min.u}
     \int_{\spt \phi} f \left( x, \delta u, D^m u \right) \, dx \le 
     \int_{\spt \phi} f \left( x, \delta u + \delta \phi, D^m u + 
     D^m \phi \right) \, dx,
\end{equation}
for all $\phi \in W^{m,1}_0 \left( \Omega ; \bR^N \right)$ with 
$\left|D^m \phi\right|^{p\left(\cdot\right)} \in L^{1}_{loc}\left(\Omega\right)$ 
and $\spt \phi \Subset \Omega$.

We will prove that under additional continuity assumptions on the functional with respect to the first and second variable and provided the exponent function $p$ is H\"older
continuous, the minimizer $u$ is regular on an open set $\Omega_0 \subset \Omega$
of full Lebesgue measure, in the sense that $D^m u$ is H\"older continuous on $\Omega_0$.

Key to the proof is an estimate for the so--calles excess
$\Phi$, which is defined by
\begin{equation*}
\Phi^2(x_0,\rho) \equiv \midint_{B(x_0,\rho)} \left| V_{p_2}(D^m u) - \left(
V_{p_2}(D^m u)\right)_{x_0,\rho} \right|^2\, dx,
\end{equation*}
with $V_p:\bR^{\bN} \to \bR^{\bN},\ z \mapsto \left(1+|z|^2\right)^{\frac{p-2}{4}}z$
and where $p_2$ is the maximal exponent on a suitable ball.
The function $\Phi$ provides an integral measure for the oscillations of $D^m u$ in
the ball $B_\rho$. A decay estimate for $\Phi$ --- which can be shown under certain initial smallness conditions on $\Phi$ --- leads to H\"older continuity of $D^m u$
via the integral characterization of H\"older continuous functions due to Campanato
(see \cite{Campanato:1965}). The excess-- decay estimate will be proved by the
'blow--up' technique: Supposing that the excess--decay estimate does not hold, but the
excess $\Phi$ is initially small, one finds a suitable blow--up sequence converging to the solution of a 'frozen' problem, which turns out to satisfy
a 'good' estimate, whereas the 'bad' estimate from the beginning carries over to the 
elements of the blow--up sequence. This finally leads to a contradiction, provided
that the convergence of the blow--up sequence is strong enough. 

The proof of this regularity theorem is one part of the authors Ph.D. thesis (see
\cite{Habermann:Diss}). It follows the ideas of Acerbi \& Mingione \cite{Acerbi:2001}, who proved the result for first order functionals of the type
\begin{equation*}
     \cF[u] := \int_\Omega f(x,Du)\, dx.
\end{equation*}
Nevertheless there are some additional difficulties to overcome due to the higher 
order case and the dependence of the integrand on lower order derivatives $D^k u\ (k 
= 0,\ldots,m-1)$. 

As mentioned above, in order to show an excess decay estimate, the minimizer of the original problem 
will be 'compared' to the solution of a corresponding 'frozen' problem. Since the
frozen problem turns out to be a system of linear PDEs with constant coefficients, 
its solution can be shown --- as in the first order case --- to be smooth and satisfy
a suitable estimate. 

However additional difficulties come up by the choice of the blow--up sequence. In
order to guarantee the convergence of a subsequence, one needs to assume the boundedness 
of alle mean values $(D^k u)_{x_0,R},\ k = 0,\ldots,m$ on the ball $B(x_0,R)$. This
leads to an explicit restriction on the regular set $\Omega_0$.

On the other hand, for proving a priori higher integrability of $D^m u$, one needs 
a suitable Caccioppoli inequality. Therefore it will be necessary to test by a function $\phi$ which is an element of the space $W^{m,p}_0$, i.e. for which the
mean values of all derivatives $D^k \phi$ up to order $m-1$ vanish. Therefore
the test function has the form $\eta^{mp} (u-P)$, where $\eta$ is a cut--off function 
and $P$ denotes the unique polynomial of order $m-1$ whose coefficients are chosen to 
satisfy
\begin{equation*}
     \midint_{B_\rho} D^k(u-P)\, dx = 0, \quad \mbox{ for }\ k=0,1,\ldots,m-1.
\end{equation*}
Existence and uniqueness of such polynomials are well known and can be found for example in \cite{Giaquinta:1979}.

The author should mention that for the sake of brevitiy, some of the proofs of the
preliminary statements are shortened very much. At many points there are only pointed out the differences to the first order case. All of the statements are proved in a
careful and extensive way in \cite{Habermann:Diss}. 

\subsection*{Acknowledgements}

I would like to thank Prof. Dr. Frank Duzaar and Prof. Dr. Giuseppe Mingione for 
encouraging me to study this kind of regularity problems and for the fruitful and intensive discussions about functionals with $p(x)$-- growth.

\section{Notations and Setting}\label{sec:setting}

We consider minimizers $u \in W^{m,1}_{loc}(\Omega;\bR^N)$ of the functional $\cF$,
where the integrand function $f$ is a Carath\'eodory function which satisfies the
growth condition
\begin{equation}\label{Wachst.2}
     L^{-1} \left( \mu^2 + | z |^2 \right) ^{\frac{p(x)}{2}} 
     \le f(x,\xi,z) \le L \left( \mu^2 + |z|^2 \right)^{\frac{p(x)}{2}},
\end{equation}
for all $x \in \Omega$, $\xi \in \bR^{\cM}$, $z \in \bR^{\cN}$, and $f$ is uniformly strictly quasiconvex with respect to the third variable, i.e. \begin{equation}\label{Ellipt.}
\begin{aligned}
     \int_Q ( f( x_0,\xi_0,z_0 + D^m \phi ) - 
          &f( x_0,\xi_0,z_0 )) \, dx,\\
     &\ge \frac{1}{L} \int_Q \left( \mu^2 + \left| z_0 \right|^2 
          + \left| D^m \phi \right|^2 \right)^{\frac{p\left(x_0\right) - 2}{2}}
          \left| D^m \phi \right|^2 \, dx,
\end{aligned}
\end{equation}
for all $z_0 \in \bR^{\cN}, x_0 \in \Omega, \xi_0 \in 
\bR^{\cM},\phi \in C_0^{\infty} 
\left( Q , \bR^N \right)$ with $0 < \mu \le 1$ and $L \ge 1$. Here $Q = \left( \left] 
0,1 \right[ \right)^n$ denotes the unit cube in $\bR^n$.

Furthermore we demand a continuity condition in the first variable of $f$ of the
type
\begin{equation}\label{Stetigk.x}
\begin{aligned}
     | &f(x,\xi,z)  - f \left(x_0,\xi,z \right) |\\
     &\le L \omega_1 
          \left( \left| x - x_0 \right| \right) \left( \left( \mu^2 + 
          |z|^2\right)^{\frac{p(x)}{2}} + \left( \mu^2 + |z|^2 
          \right)^{\frac{p\left(x_0\right)}{2}}\right)
          \cdot \left( 1 + \left| \log \left( 
          \mu^2 + |z|^2 \right)\right| \right),
\end{aligned}
\end{equation}
for all $x, x_0 \in \Omega, \xi \in \bR^{\cM},
z \in \bR^{\cN}$. 
$\omega_1:\left( 0, \infty \right) \to \left( 0, \infty \right)$ denotes the modulus 
of continuity of the function $p$, i.e. $\omega_1$ is non decreasing, concave, 
continuous and
\begin{equation}\label{Stet.Mod.omega1a}
     \lim \limits_{R \downarrow 0} \omega_1 \left( R \right) = 0,
\end{equation}
as well as
\begin{equation}\label{Stet.Mod.omega1b}
     \left| p(x) - p(y) \right| \le \omega_1 \left( \left| x-y \right| \right).
\end{equation}

Finally we suppose the function $f$ to be continuous in the second variable, i.e.
\begin{equation}\label{Stetigk.xi}
     \left|f\left(x,\xi,z\right) - f\left(x,\xi_0,z\right)\right| \le L \omega_2
     \left(\left|\xi - \xi_0\right|\right)\left(\mu^2 + \left|z\right|^2
     \right)^{\frac{p(x)}{2}}
\end{equation}
for all $\xi,\xi_0 \in \bR^{\cM}, x \in 
\Omega, z \in \bR^{\cN}$.
Without loss of generality we assume that $\omega_2$ is concave, bounded and therefore
subadditive.

The main statement is the following


\begin{theorem}\label{Hauptsatz} 
Let $u \in W^{m,1}_{loc}\left(\Omega;\bR^N\right)$ be a local minimizer of the 
functional $\cF$, where $f$ denotes a function of the class $C^2$ with respect to the 
variable $z$, which satisfies the growth, quasiconvexity and continuity assumptions (\ref{Wachst.2}) to (\ref{Stetigk.xi}). Furthermore let $\alpha \in (0,1]$ and the
moduli of continuity $\omega_1$ and $\omega_2$ satisfy
\begin{equation}\label{Stet.Mod.stark}
     \omega_1(R) + \omega_2\left(R\right) \le L R^{\alpha},
\end{equation}
for all radii $R \le 1$. Then there exists an open subset 
$\Omega_0 \subset \Omega$ of full measure, i.e. ${\mathcal 
L}^n\left(\Omega \setminus \Omega_0\right) = 0$, such that $D^m u$ is locally H\"older continuous in $\Omega_0$ for some H\"older exponent.
\strut\hfill $\blacksquare$
\end{theorem}

\begin{remark}\label{Wachst.Df}
The regularity result stated in Theorem \ref{Hauptsatz} does
not assume any further growth assumptions, especially not on the second derivatives
of the integrand function $f$.

If (\ref{Wachst.2}) holds, condition (\ref{Ellipt.}) implies the following growth 
condition for the first derivative $Df$ (if it exists):
\begin{equation}\label{Wachst. Df 2}
     \left|Df\left(x_0,\xi_0,z\right)\right| \le c\left(1+\left|z\right|
     \right)^{p\left(x_0\right) - 1},
\end{equation}
with a constant $c \equiv c\left(L,p(x_0)\right)$, for all $z \in 
\bR^{\cN}$ and $x_0 \in \Omega, \xi_0 \in \bR^{\cM}$. Additionally if $1 < \gamma_1 \le p(x_0) \le \gamma_2 < + \infty$,
then the constant $c$ depends only on $n$ and $\gamma_2$.
\hfill $\blacksquare$
\end{remark}

\subsubsection*{Remarks on the notation}

In the whole paper
$\Omega \subset \bR^n, (n\ge 2)$ denotes a bounded domain in the space $
\bR^n$ and $B\left(x,R\right) \equiv B_R(x)$ the open ball $\left\{ y \in 
\bR^n:\ \left|x-y\right| < R\right\}$.
The Lebesgue measure of a measurable set $A$ is abbreviated by $|A| \equiv {\mathcal L}^n(
A)$. For a locally integrable function $u \in L^1_{loc}(\Omega)$ we define the mean value on the ball $B$ by
\begin{equation*}
     \left(u\right)_{x_0,R} := \midint_{B\left(x_0,R\right)} 
        u(x)\, dx = \frac{1}{|B(x_0,R)|}\int_{B\left(x_0,R\right)} 
        u(x)\, dx.
\end{equation*}
In the case the centre of the ball is obvious from the context, we will often
just write $B_R$ or $B$ instead of $B(x_0,R)$, $(u)_R$ instead of $(u)_{x_0,R}$
respectively.

The letter $c$ denotes a constant which will not necessarily be the same at different places
in the work and which may sometimes change from line to line. Constants that will
be referred to at other points of the work, will be signed in a unique way, mostly by
different indices. In the case we want to emphasise the fact that a constant changes
from one line to another, we will label this by mathematical accents, as for example
$\tilde{c}$ or $\bar{c}$. For the survey we will not specify the dependencies of the 
constants in between the estimates, but of course at the end of them.

For $\Omega \subset \bR^n, p > 1$, let $L^{p}(\Omega;\bR^N)$ be the well known Lebesgue
space to the power $p$. For $m \in \bN$ we define the Sobolev space 
\begin{equation*}
     W^{m,p}\left( \Omega; \bR^N \right) := \left\{ u \in L^p\left( \Omega, 
     \bR^N \right): D^{\alpha}u \in L^p\left( \Omega\right) \mbox{ for } 0 
     \le \left|\alpha\right|\le m \right\},
\end{equation*}
with the multi-index $\alpha = (\alpha_1,\ldots,\alpha_n) \in \bN^{n}$ and the
abbreviations $|\alpha| := \alpha_1 + \ldots + \alpha_n$ and $D^{\alpha} u :=
D^{\alpha_1}_1\ldots D^{\alpha_n}_n u$. Furthermore let $W^{m,p}_0(\Omega;\bR^N)$
denote the closure of $C^{\infty}(\Omega;\bR^N)$ in the space $W^{m,p}(\Omega;\bR^N)$.

\section{Basic tools and preliminary results}

\subsection{General assumptions}

For the whole paper we will assume that
\begin{equation}\label{Stet.Mod.3}
     \limsup \limits_{R \downarrow 0} \omega_1 \left( R \right) 
     \log \frac{1}{R} < \infty,
\end{equation}
which is a weaker condition than condition (\ref{Stet.Mod.stark}). Therefore without 
loss of generality we can also assume, that for all $R \in (0,1]$ there 
holds
\begin{equation}\label{Stet.Mod.4}
     \omega_1 \left( R \right) \log \frac{1}{R} \le L,
\end{equation}
where $L$ is the constant from the growth condition (\ref{Wachst.2}). As all our 
results are local, we can furthermore assume that
\begin{eqnarray}
     1 < \gamma_1 \le p(x) \le \gamma_2 & < & \infty \ \ \ \ \mbox{ for all }\ x \in 
         \Omega, \label{Beschr.p}\\
     \int_{\Omega} \left| D^m u \right|^{p(x)} \, dx & < & \infty.
         \label{Integr.u}
\end{eqnarray}

\subsection{Higher integrability results}

\subsubsection{\bf Higher integrability of \boldmath$|D^m u|^{p(\cdot)}$\unboldmath}

First we prove higher integrability for $|D^m u|^{p(\cdot)}$. A similar result for
first order functionals ($m=1$) was shown in \cite{Acerbi:2001}.

\begin{Lem}\label{High.Int.p(x)}
Let $O \subset \Omega$ be open and $u \in W_{loc}^{m,1} \left(O, \bR^N \right)$ 
a local minimizer of the functional
\begin{equation}
     w \mapsto \int_{O} f \left(x, \delta w(x), D^m w(x) \right) \, dx,
\end{equation}
where $f: O \times \bR^{\cM} \times 
\bR^{\cN} \to \bR$ satisfies the growth and continuity assumptions
(\ref{Wachst.2}) and (\ref{Stetigk.x}).
Furthermore assume that
\begin{equation}\label{p(x)-Integr.O}
     \int_{O} \left| D^m u \right|^{p(x)} \, dx \le M < \infty.
\end{equation}
Then there exist $\delta,c \equiv \delta,c\left( n, \gamma_1, \gamma_2, 
L, M, m \right) > 0$ and a radius $R_0 \equiv R_0 ( n, \gamma_1, 
\omega(\cdot))$, such that for every ball $B_R \Subset O$ with $R \le R_0$ there holds:
\begin{equation}\label{High.Int.p(x).2}
     \Biggl( \midint_{B_{\frac{R}{2}}} \left| D^m u \right|^{p(x)
     \left( 1 + \delta \right)} \, dx \Biggr)^{\frac{1}{1+\delta}} \le c 
     \Biggl( \midint_{B_R} \left| D^m u \right|^{p(x)} \, dx + 1 
     \Biggr).
\end{equation}
\hfill$\blacksquare$
\end{Lem}

\begin{proof}
The proof of this result is more or less standard. Therefore we will show only
the main steps, especially pointing out the additional difficulties in the higher
order case. 

 Let $t,s \in (0,1)$ be such that $\frac{R}{2} <t <s < R \le 1$ and moreover
\begin{equation}
     p_1 := \min\left\{p(x): x \in B_{4R} \right\},\ \ \ \ \ p_2 := \max\left\{
     p(x): x \in B_{4R}\right\}.
\end{equation}
We choose a cut-off function $\eta \in C_c^{\infty} \left( B_R \right)$ with the
following properties:
\begin{equation}\label{Abschn.fkt.1}
     0 \le \eta \le 1 \ \ \mbox{ on }\ B_R,\ \ \ \ \eta \equiv 0 \ \ \mbox{ 
     outside}\ B_s,\ \ \ \ \eta \equiv 1 \ \ \mbox{ on }\ B_t,
\end{equation}
\begin{equation}\label{Abschn.fkt.2}
     \left| D^k \eta \right| \le \frac{c_{co}(m)}{\left( s-t \right)^k} \ \ \ \mbox{ 
     for all }\ \ k = 1, \ldots , m.
\end{equation}
Furthermore we set 
\begin{equation}\label{Testfkt.1}
     \varphi(x) := \eta(x) \left( u(x) - P(x) \right),
\end{equation}
where $P: \Omega \to \bR^N$ denotes the unique polynomial of degree $m-1$ which satisfies
\begin{equation}\label{Eigensch.l}
     \left( D^k\left(u-P\right)\right)_R \equiv \midint_{B_R} D^k \left( 
     u-P \right) \, dx = 0 \ \ \ \ \mbox{ for } \ \ k = 0,\ldots , m-1.
\end{equation}
Finally we define 
\begin{equation}\label{Testfkt.2}
     v := u - \varphi = u - \eta \left( u-P \right) = \left( 1 - \eta \right) u 
     + \eta P.
\end{equation}
Then we have $v \equiv u$ on $B_R \setminus B_s$, $v \equiv P$ on $B_t$. By the
growth condition (\ref{Wachst.2}), minimality of $u$ and the definition of 
$v$ we easily see
\begin{equation*}
     \int_{B_t} \mskip-5mu \left| D^m u \right|^{p(x)} \, dx
     \le L^2 \int_{B_s} \mskip-5mu \left( \left| D^m \left( u - \eta 
     \left( u-P \right) \right)\right|^{p(x)}+1\right) \, dx.
\end{equation*}
We further estimate the right hand side, using the definition of the cuf--off function
and the global bounds \eqref{Beschr.p} of the exponent function, which yields
\begin{equation*}
     \int_{B_t} \left| D^m u \right|^{p(x)} \, dx \le c\left[\int_{B_s 
          \setminus B_t} \left|D^m u\right|^{p(x)} \, dx
          + \sum\limits_{k=1}^m\int_{B_s}\left| \frac{D^{m-k}\left(u-P\right)}{\left( t-s\right)^k}
          \right|^{p_2}\, dx + |B_R|\right],
\end{equation*}
with a constant $c \equiv c(m,\gamma_2,L)$.
With the standard 'hole filling' technique and application of Lemma \ref{Techn.Lemma1} we end up with
\begin{equation*}
\int_{B_{\frac{R}{2}}} \left|D^m u \right|^{p(x)}\, dx \le c
     \Biggl[\sum\limits_{k=1}^m \int_{B_R}\left|
     \frac{D^{m-k}\left(u-P\right)}{R^k}\right|^{p_2}
     \, dx + |B_R| \Biggr].
\end{equation*}
By the choice of our polynomials $P$ we are in a position to apply Poincar\'e's
inequality iteratively to the integrals on the right hand side, which finally leads us to the following Caccioppoli--type inequality:
\begin{equation}
 \midint_{B_{\frac{R}{2}}}\left|D^m u \right|^{p(x)} \, dx \le c_3\left[ \midint_{
     B_R}\left|\frac{D^{m-1}\left(u-P\right)}{R}\right|^{p_2}\, dx + 1\right],
\end{equation}
with a constant $c \equiv c(n,m,L,\gamma_2)$.
Now we proceed exactly in the same way as shown in \cite{Acerbi:2001}, first 
applying Sobolev-Poincar\'e's inequality to the right hand side, then choosing
$\vartheta := \min \{\sqrt{\frac{n+1}{n}}, \gamma_1\}$ and 'localizing', i.e.
choosing the radius $R$ so small that $\omega(8R) \le \vartheta - 1$, which allows
us to pass over from $p_2$ to the variable exponent $p(x)$. In conclusion we 
end up with a reverse H\"older inequality of the following type:
\begin{equation}\label{Rev.H}
\midint_{B_{\frac{R}{2}}} \left|D^m u\right|^{p(x)} \, dx \le c\left[ 
     \left(\ \midint_{B_R} \left|D^m u\right|^{\frac{p(x)}{\vartheta}}\, dx
     \right)^{\vartheta} + 1\right],
\end{equation}
for all $B_R \Subset \Omega$ with $8R \le 1,\ \omega\left( 8R\right) 
\le \vartheta -1$, where $\vartheta \equiv \min \left\{ \gamma_1, \sqrt{\frac{n+1}{
n}} \right\}$, $c \equiv c\left(n,\gamma_1,\gamma_2,L,M,m\right)$.

Now let $R_0$ be the maximal $R$ with these properties, thus $8R_0 \le 1,\ \omega
\left( 8R_0\right) \le \vartheta -1$, i.e. $R_0 = R_0 \left( n,\gamma_1,
\omega(\cdot)\right)$. Then (\ref{Rev.H}) holds for any $B_R \Subset \Omega$ with $R \le 
R_0$. The statement of the lemma follows by an application of Gehring's lemma
(see \cite[Chapter V]{Giaquinta:1983}) with $f \equiv \left| D^m u\right|^{ 
p(x)}$.
\end{proof}

\subsubsection{\bf An up--to--the--boundary higher integrability result}

The following result up to the boundary is stated in \cite{Acerbi:2001} and we
will follow the ideas from there, whereas we should mention that the statement which
is proved here is slightely better, since we end up with the same radius of the balls
on both sides of the inequality. This is due to a global version of the Gehring
lemma, stated in \cite{Duzaar:2004}.

\begin{Lem}\label{High.Int.Rand}
Let $\Omega \subset \bR^n$ be an open set, $p$ a constant exponent, $1\le\gamma_1 \le 
p \le \gamma_2$, $R>0,\ B_R \subset \Omega$. Further let $g:\Omega \times 
\bR^{\cN} \to \bR$ be continuous such that
\begin{equation}\label{Wachst.g}
     L^{-1}|z|^p \le g(x,z) \le L\left(|z|^p + a(x)\right),
\end{equation}
for all $z \in \bR^{\cN}$ with $L \ge 1$, $0 < a \in L^{\gamma}
\left(B_R\right)$ for some $\gamma > 1$. For given $h \in W^{m,q}\left(B_R\right)$ 
with $q>p$,
let $v$ be the solution of the minimization problem
\begin{equation}\label{Var.Int}
     \min\left\{\int_{B_R} g\left( x,D^m w\right)\, dx,\ \ \ w \in 
     h + W_0^{m,p}\left(B_R\right) \right\}.
\end{equation} 
Then there exists $\varepsilon = \varepsilon\left(n,\gamma_1, \gamma_2, L,m\right) 
\in \left]0,\tilde{m}\right[$ with $\tilde{m} = \min\left\{\gamma - 1,
\frac{q}{p}-1\right\}$ and a constant $c$ depending only on $n,\gamma_1,\gamma_2,L,m$, such that
\begin{equation}
\begin{aligned}
     \midint_{B_R}|D^m v|^{p\left(1+\varepsilon
          \right)} \, dx
     \le c \Biggl[&\Biggl(\midint_{B_R} \left|D^m v\right|^p \, dx 
          \Biggr)^{1+\varepsilon}\\
     &+ \Biggl(\midint_{B_R}\left|D^m h 
          \right|^{p\left(1+\tilde{m}\right)}\, dx\Biggr)^{\frac{1+\varepsilon}{ 
          1+\tilde{m}}}
      + \Biggl(\midint_{B_R}a^{1+\tilde{m}}\, dx
          \Biggr)^{\frac{1+\varepsilon}{1+\tilde{m}}}\Biggr].
\end{aligned}
\end{equation}
\hfill$\blacksquare$
\end{Lem}

\begin{proof}
We distinguish the interior situation and the situation on the boundary. In the
interior situation $B_\rho = B_\rho(x_0) \subset B_R$ the proof is similar to the 
proof of Lemma \ref{High.Int.p(x)} . It is even more simple since the exponent $p$ is constant. 
Proceeding in the standard way and additionally using Poincar\'e's inequality in the
same way as in the proof of Lemma \ref{High.Int.p(x)} , we end up with the following Caccioppoli--
inequality:
\begin{equation*}
\midint_{B_{\frac{\rho}{2}}} \left|D^m v\right|^p\, dx \le c \left[\midint_{B_{
     \rho}}\left|\frac{D^{m-1}\left(v-P\right)}{\rho}\right|^p\, dx + 
     \midint_{B_{\rho}}a\, dx \right],
\end{equation*}
with $c \equiv c(n,L,m,\gamma_2)$.

In the boundary situation $B_\rho = B_\rho(x_0)$ with $x_0 \in \partial B_R$ we proceed in the following way: Let $\frac{\rho}{2} < t < s < \rho$ and $\eta$ a
cut--off function as in the proof of Lemma \ref{High.Int.p(x)} (of course we have to take $\rho$
instead of $R$ here). We define 
\begin{equation*}
     w := v - \eta (v-h).
\end{equation*}
Then there holds $w = v$ on $\partial B_\rho$ and on $B_\rho \setminus B_s$ and
$w = h$ on $B_t$. Furthermore by (\ref{Var.Int}) we have $D^kv = D^kh$ on $\partial B_R$ 
for $k=0,\ldots,m-1$. Therefore by minimality of $v$, the imposed boundary condition, the imposed growth condition and the special form of the
function $w$ we obtain, after again applying the 'hole filling' technique and using again Lemma \ref{Techn.Lemma1}, the estimate (Note that we use the 
notation $B_r^+ := B_r \cap B_R$)
\begin{equation*}
     \int_{B_{\frac{\rho}{2}}^+}\left|D^m v\right|^p\, dx \le c
     \Biggl[\sum\limits_{k=1}^m \int_{B_{\rho}^+}\left|\frac{D^{m-k}
     \left(v-h \right)}{\rho^k}\right|^p\, dx + \int_{B_\rho^+} |D^m h|^p + 
     \int_{B_\rho^+} a\, dx \Biggr],
\end{equation*}
with $c \equiv c(m,L,\gamma_2)$. Now we define the function
\begin{equation*}
\tilde{v} := \left\{ \begin{array}{ll} D^{m-k}(v-h) & \mbox{ on } B_{\rho}^+,\\
     0 & \mbox{ on } B_{\rho} \setminus B_{\rho}^+ .\\ \end{array}\right. 
\end{equation*}
Since $v-h \in W^{m,p}_0(B_R)$, we can iteratively apply Poincar\'e's inequality in 
the version of \cite{Ziemer:1989}, Corollary 4.5.3, to $\tilde{v}$ in combination
with H\"older's inequality to conclude with the following Caccioppoli--inequality at
the boundary:
\begin{equation}\label{Cacc.,Rand}
     \midint_{B_{\frac{\rho}{2}}^+}\left|D^m v\right|^p\, dx 
     \le c \Biggl[\midint_{B_{\rho}^+}\left|\frac{D^{m-1}\left(
          v-h\right)}{\rho}\right|^p\, dx + \midint_{B_{\rho}^+}
          \left|D^m h\right|^p\, dx
     + \midint_{B_{\rho}^+}a\, dx \Biggr].
\end{equation}

To conclude inequalities of reverse H\"older type, we estimate in the interior
situation via Sobolev-Poincar\'e's and H\"older's inequalities in the standard
way. In the situation at the boundary we define
\begin{equation*}
     \tilde{v}:=\left\{\begin{array}{lcl} D^{m-1}\left(v-h\right) & \mbox{ on }& B_{
     \rho}^+,\\0 & \mbox{ on }& B_{\rho}^- := B_{\rho} \setminus B_{\rho}^+ ,\\ 
     \end{array} \right.
\end{equation*}
and apply the Sobolev-Poincar\'e inequality in the version of 
\cite[Corollary 4.5.3, p. 452]{Ziemer:1989} to obtain (note that $|B_\rho^-| \ge
1/2 |B_\rho|$)
\begin{eqnarray*}
     \int_{B_{\rho}^+} \left| D^{m-1}\left(v-h\right)\right|^p\, dx
     & = &\int_{B_{\rho}}\left|\tilde{v}\right|^p\, dx
     \le c\frac{|B_{\rho}|}{|B_{\rho}^-|}\Biggl(\int_{B_{\rho}} 
          |D\tilde{v}|^{\frac{np}{n+p}}\, dx\Biggr)^{\frac{n+p}{n}}\\
     &\le&c\Biggl(\int_{B_{\rho}^+}\left|D^{m}\left(v-h\right)\right|^{
          \frac{np}{n+p}}\, dx\Biggr)^{\frac{n+p}{n}},
\end{eqnarray*}
with $c \equiv c(n,\gamma_2)$. Substituting this into \eqref{Cacc.,Rand} and 
subsequently applying H\"older's inequality we finally end up with the 
following reverse H\"older inequalities:
\begin{equation}\label{HI.bdry.RH_Rand}
     \midint_{B_{\frac{\rho}{2}}^+}\left|D^m v\right|^p\, dx \le \tilde{c} \Biggl[
     \Biggl(\midint_{B_{\rho}^+}\left|D^m v\right|^{p\chi}\, dx\Biggr)^{1/\chi}+ 
     \midint_{B_{\rho}^+}\left( \left|D^m h\right|^p + a\right)\, dx\Biggr],
\end{equation}
and 
\begin{equation}\label{HI.bdry.RH_innen}
     \midint_{B_{\frac{\rho}{2}}} \left|D^m v\right|^p\, dx \le \tilde{c}
          \Biggl[\Biggl(\midint_{B_{\rho}}\left|D^m v\right|^{p\chi}\, dx\Biggr)^{
          1/\chi} + \midint_{B_{\rho}}a\, dx \Biggr],
\end{equation}
with $\tilde{c} \equiv \tilde{c}(n,m,L,\gamma_2)$ and $\chi \equiv \frac{n}{n+p} < 1$. Note that \eqref{HI.bdry.RH_Rand} holds for $\rho \le R$ and \eqref{HI.bdry.RH_innen} for
all $B_{\rho} \subset B_R$. Therefore we can apply
the global version of the Gehring lemma in \cite[Theorem 2.4]{Duzaar:2004}, with the functions
\begin{eqnarray*}
     g & := & \left|D^m v\right|^{p\chi},\\
     f & := & \left(\left|D^m h\right|^p + a \right)^{\chi}
\end{eqnarray*}
This guarantees the existence of $\varepsilon \equiv \varepsilon\left(n,p,m,L,
k_{\Omega}\right) \in \left]0,\tilde{m}\right[$ with $\tilde{m} := 
\min\left\{\gamma - 1, \linebreak \frac{q}{p} - 1 \right\}$ such that
\begin{equation*}
     \Biggl(\midint_{B_R}\left|D^m v\right|^{p\chi\tilde{q}}\, dx\Biggr)^{\frac{1}{
          \tilde{q}}}
     \le c\Biggl[\Biggl(\midint_{B_R}\left|D^m v\right|^p\, dx
          \Biggr)^{\chi} + \Biggl(\midint_{B_R}\left(\left|D^m h\right|^p + a
          \right)^{\chi\tilde{q}}\, dx \Biggr)^{\frac{1}{\tilde{q}}}\Biggr],
\end{equation*}
for all $\tilde{q} \in \left[\frac{1}{\chi}, \frac{1}{\chi}\left(1+\varepsilon
\right)\right]$, with $c \equiv c(n,m,L,p)$. Choosing $\tilde{q} := \frac{1}{\chi}
\left( 1 + \varepsilon \right)$ and raising the resulting inequality to the power
$\left(\chi p\right)^{-1}$ yields the desired estimate.

\begin{remark}
As one can easily deduce from the proof of Theorem 2.4 in \cite{Duzaar:2004}, 
the constant in the estimate above can be replaced by a constant depending only on $n, m, L, \gamma_1$ and $\gamma_2$.
\end{remark}
\end{proof}

\subsubsection{\bf Higher integrability for an almost minimizer}
The following lemma will be needed to have higher integrability for the comparison
function in the blow--up procedure, obtained by Ekeland's principle. The result in
first order case can be found in \cite{Carozza:1998} for the case $1 < p < 2$ and
in \cite{Acerbi:1987} for the case $p \ge 2$. Since the proof in the higher order
case is only a slight modification of the proofs in case of first order, we do not
give it here and refer the reader to \cite{Habermann:Diss} for a detailed discussion.

\begin{Lem}\label{High.Int.3}
Let $p$ be constant, $1 < \gamma_1 \le p \le \gamma_2$, $\tilde{L} > 1,\ 0 < 
\lambda < 1$ and $g:\bR^{\cN} \to \bR$ be a continuous function 
satisfying the following conditions:
\begin{eqnarray}
     \left|g(z)\right| &\le& \tilde{L} \left( 1 + \lambda^2\left|z\right|^2
          \right)^{\frac{p-2}{2}}\left|z\right|^2 = \tilde{L} \lambda^{-2} 
          \left|V_p(\lambda z)\right|^2,\\
     \int_{B_1} g\left(D^m\phi\right)\, dx &\ge& \tilde{L}^{-1}\int_{B_1} 
          \lambda^{-2} \left|V_p\left(\lambda D^m \phi\right)
          \right|^2\, dx,
\end{eqnarray}
for all $\phi \in W_0^{m,p}\left(B_1;\bR^N\right)$. Moreover let $0 < \mu < 1$ and 
$\bar{u} \in W^{m,p}\left(B_1;\bR^N\right)$ such that there holds
\begin{equation}
     \int_{B_1}g\left(D^m \bar{u}\right)\, dx \le \int_{B_1}\left(g\left(D^m
     \bar{u} + D^m \phi\right) + \mu\left|D^m \phi\right|\right)\, dx,
\end{equation}
for all $\phi \in W_0^{m,p}\left(B_1;\bR^N\right)$. Then there exist constants
$c,\delta_2 \equiv c,\delta_2 (n,N,m,\gamma_1, \gamma_2,\tilde{L})$
independent of $R,\lambda,\bar{u},g, \mu$ such that for every ball
$B_{3R} \subset B_1$ there holds
\begin{equation}
     \Biggl(\midint_{B_{\frac{R}{2}}} \left|V_p\left(\lambda D^m \bar{u}
     \right)\right|^{2\left(1+\delta_2\right)}\, dx\Biggr)^{\frac{1}{1+ 
     \delta_2}}\le c \midint_{B_{3R}}\left(\lambda^2 \mu + \left|V_p\left( 
     \lambda D^m \bar{u}\right)\right|^2\right)\, dx.
\end{equation}
\hfill$\blacksquare$
\end{Lem}

\subsubsection{\bf Technical Lemma}

The following technical lemma, concerning the Taylor approximation of the function
$f$ in the point $x_0$ will be needed later in the proof of the main theorem. 
For the proof we refer the reader to \cite{Habermann:Diss}.
\begin{Lem}\label{Techn.Lemma2}
Let $M>1$. For $x_0 \in \Omega,\ U = (U_1,\ldots,U_m) \in 
\bR^{\cM}, A \in \bR^{\cN}$ 
with $\left|A\right|\le M, |U_i| \le M\ (i=1,\ldots,m)$ and $\lambda > 0$ let
\begin{equation*}
     f_{U,A,\lambda}\left(\zeta\right) := \lambda^{-2}\left[f\left(x_0,U,A + 
     \lambda \zeta\right) - f\left(x_0,U,A\right) - \lambda Df\left(x_0,U,A
     \right) \zeta \right],
\end{equation*}
where the function $f$ satisfies \eqref{Wachst.2} and \eqref{Ellipt.}.
Then there exists a constant $\bar{L} \equiv \bar{L}\left(\gamma_1,\gamma_2,
L,M\right)$ such that
\begin{eqnarray*}
     \mbox{a)} && f_{U,A,\lambda}\left(\zeta\right) \le \bar{L}\left(1+\lambda^2
           \left|\zeta\right|^2\right)^{\frac{p\left(x_0\right)-2}{2}}\left|
          \zeta\right|^2 = \bar{L} \lambda^{-2} \left|V_{p\left(x_0\right)}
          \left(\lambda \zeta\right)\right|^2,\\
     \mbox{b)} && \int_{B_1} f_{U,A,\lambda}\left(D^m \phi\right)\, dx 
          \ge \bar{L}^{-1} \int_{B_1} \lambda^{-2}\left|V_{p\left(x_0
          \right)}\left(\lambda D^m \phi\right)\right|^2 \, dx,\\
\end{eqnarray*}
for all $\phi \in W^{m,p\left(x_0\right)}_0\left(B_1;\bR^N\right)$.
\hfill$\blacksquare$

\end{Lem}

\subsection{Ekeland variational principle}

The following lemma will play a central role in the proof of the main theorem. In
the blow--up procedure we will need this variational principle to compare the
minimizer $u$ to an almost minimizer of a suitable 'frozen' problem. The lemma
is cited from \cite{Ekeland:1979}.

\begin{Lem} \label{eke_orig}
Let $(\cX,d)$ be a complete metric space and $\cG
: \cX \to (-\infty,+\infty]$ a lower
semicontinuous functional such that $\inf_{\cX}
\cG$ is finite. Given $\eps>0$ let $u\in \cX$ be
such that $\cG(u) \le \inf_{\cX} \cG +\eps$. Then there
exists $w \in \cX$ such that
\begin{eqnarray*}
d(w,u) &\leq& 1, \\
\cG(w) &\leq& \cG(u),\\
\cG(w) &\leq& \cG(v) +\eps d(v,w), \quad
\text{for any $v\in \cX$.}
\end{eqnarray*}
\end{Lem}

\subsection{The function $V_p$}

Let the function $V \equiv V_p:\bR^l \to \bR^k$ be defined by
\begin{equation}\label{Def.V_p}
     V_p(z) = \left(1+|z|^2\right)^{\frac{p-2}{4}}z.
\end{equation}
We recall algebraic properties of the function
$V_p$ (for a proof of the properties see
e.g. \cite{Carozza:1998}).

\begin{Lem}\label{prop.V}
Let $p>1$ and let $V \equiv V_p: \bR^k \to \bR^k$
be as in \eqref{Def.V_p}. Then for any $z, \eta\in \bR^k$

\begin{enumerate}[i)]
\item $|V(tz)| \leq \max\{t,t^{p/2}\}|V(z)|$,\quad for any  $t>0$;
\item $|V(z+\eta)| \leq c\Big(|V(z)| +|V(\eta)|\Big)$;
\item
$$
c^{-1}|z-\eta| \leq
\frac{|V(z)-V(\eta)|}{(1+|z|^2+|\eta|^2)^{(p-2)/4}}
\leq c|z-\eta|;
$$
\noindent Moreover for any $z\in \bR^k$
\item
\begin{align*}
\text{if $p\in(1, 2)$:} &\quad \frac{1}{\sqrt{2}}\min\{|z|,|z|^{p/2}\} \leq |V(z)| \leq
\min\{|z|,|z|^{p/2}\}, \\
\text{if $p\geq 2$:} &\quad \max\{|z|,|z|^{p/2}\} \leq |V(z)| \leq
\sqrt{2}\max\{|z|,|z|^{p/2}\},
\end{align*}
\item
\begin{align*}
\text{if  $p \in(1,2)$:}& \quad |V(z)-V(\eta)| \leq c |V(z-\eta)|, \quad
\text{for any $\eta \in \bR^k$} \\
\text{if  $p \geq 2$:}& \quad |V(z)-V(\eta)| \leq c(M) |V(z-\eta)|, \quad
\text{for $|\eta| \leq M$}
\end{align*}
\item
\begin{align*}
\text{if  $p \in(1,2)$:}& \quad |V(z-\eta)| \leq c(M)|V(z)-V(\eta)|, \quad
\text{for $|\eta| \leq M$} \\
\text{if  $p \geq 2$:}& \quad |V(z-\eta)| \leq c|V(z)-V(\eta)|, \quad
\text{for any $\eta \in \bR^k$}
\end{align*}
\end{enumerate}
with $c(M), c\equiv c(k,p) >0$. If $1< \gamma_1 \leq p \leq
\gamma_2$ all the constants $c(k,p)$ may be replaced by a single constant
$c\equiv c(k,\gamma_1, \gamma_2)$.
\end{Lem}

\subsection{Technical Lemma II}

We now formulate a technical lemma wich is --- in a little more particular version
--- shown in \cite{Kronz:2003}. We will need this lemma in several points of the proof
of the main theorem, especially for example when proving Caccioppoli type inequalities 
by the 'hole filling' technique. For the proof in this general situation we again
refer the reader to \cite{Habermann:Diss}.

\begin{Lem}\label{Techn.Lemma1}
Let $p \ge 1$ be constant, $K \in \bN$, ${\mathcal V}_p: \bR^K \to \bR^K$ any function 
for which holds
\begin{eqnarray}
     |{\mathcal V}_p (y+z)| 
     &\le& c \left( |{\mathcal V}_p(y)| + |{\mathcal V}_p(z)|\right) \mskip+20mu 
           \mbox{ for all }\ y,z \in \bR^K \label{Techn.Lem1.1}\\
     |{\mathcal V}_p (tz)| 
     &\le& \max \{t,t^{p/2}\} |{\mathcal V}_p(z)| \mskip+25mu \mbox{ for all }\ z \in 
           \bR^K, t \in \bR \label{Techn.Lem1.2}
\end{eqnarray}
with a constant $c \equiv c(p,K) > 0$. Moreover let $0 < \theta
< 1, A_k \ge 0, a_k > 0$ for $k=0,\ldots l$, $B \ge 0$ and $f \ge 0$ be a bounded 
function satisfying
\begin{equation}\label{TL1.Ass}
     f(t) \le \theta f(s) + \sum\limits_{k=0}^l A_k \int_{B_R} 
     \left| {\mathcal V}_p\left(\frac{h_k(x)}{\left( s-t \right)^{a_k}}\right) 
     \right|^2\, dx + B
\end{equation}
for all $r < t < s < R$, where $h_k \in L^p\left(B_R\right)$ for all 
$k = 0,\ldots,l$. Then there exists a constant $c \equiv c\left(p,\theta, a_0,
\ldots,a_l\right)$ such that
\begin{equation}
     f(r) \le c \sum\limits_{k=0}^l A_k \int_{B_R} \left|{\mathcal V}_p\left(
     \frac{h_k(x)}{\left(R-r\right)^{a_k}}\right)\right|^2\, dx + B
\end{equation}
\hfill$\blacksquare$
\end{Lem} 

\subsection{Lower order derivatives}

The following lemma will enable us to estimate the $L^p$- distance between the 
term $\delta u(x)$, consisting of the derivatives of $u$ up to order $m-1$ and
the mean value $(\delta u)_{\rho}$ on the ball $B_{\rho}$ by the $W^{m,p}$-norm of $u$ .
This will be useful in several points of the proofs, utilizing the boundedness of
the mean values $(D^k u)_{\rho}$ for $k=0,\ldots m-1$.


\begin{Lem}\label{Niedr.Abl.}
Let $p \ge 1$.For $\rho > 0$, $x \in B_{\rho}(x_0)$, $u \in W^{m,p}\left(B_{\rho} 
\left(x_0 \right),\bR^N\right)$ and
\begin{equation*}
\delta u(x):= \left(u(x),Du(x),\ldots,D^{m-1}u(x)\right) \in \bR^N \times
\bR^{nN}\times \ldots \times \bR^{N\binom{n+m-2}{m-1}}
\end{equation*}
the vector whose $j$ th component contains the $j$ th weak derivative of the function $u$ in the point $x \in B_{\rho}$. Then there holds
\begin{equation*}
     \midint_{B_{\rho}} \left|\delta u(x) - \left(\delta u\right)_{\rho}\right|^p
          \, dx \le  c\sum_{k=0}^{m-1}\left[\rho^{p(m-k)}\left(\left|D^m u
          \right|^p \right)_{\rho} + \sum_{l=k+1}^{m-1}\rho^{p(l-k)}\left|\left(
          D^l u \right)_{\rho}\right|^p\right],
\end{equation*}
with a constant $c \equiv c\left(n,m,p\right)$.
\hfill$\blacksquare$
\end{Lem}

\begin{proof}
The proof consists of a sequence of elementary estimates which we will sketch in the sequel. 
Let $P_m:\bR^n \to \bR^N$ be the unique polynomial of order $m-1$, satisfying
\begin{equation}\label{LO_Poly}
     \midint_{B_{\rho}} D^k\left(u(x) - P_m(x)\right)\, dx = 0\ \ \ \mbox{ for }
     \ \ k = 0,\ldots,m-1.
\end{equation}
Then we have
\begin{eqnarray*}
     \midint_{B_{\rho}} \Bigl|D^ku - \left(D^k u\right)_{\rho}\Bigr|^p\, dx
     &\le& 2^{p-1} \left[ \midint_{B_{\rho}}\left|D^k \left(u - P_m 
          \right)\right|^p\, dx + \midint_{B_{\rho}}\left|D^kP_m - \left(D^k u 
          \right)_{\rho}\right|^p\, dx\right]\\
     & = & 2^{p-1} \left[ I_1 + I_2\right].
\end{eqnarray*}
Poincar\'e's inequality, applied $m-k$ times to $I_1$ (note \eqref{LO_Poly}) 
leads us to
\begin{equation*}
     I_1 \le c\rho^{p\left(m-k\right)}\midint_{B_{\rho}}\left|D^m u
          \right|^p\, dx.
\end{equation*}
To estimate $I_2$, we apply the explicit representation formula for the polynomial 
$P_m$ from \cite{Kronz:2005}:
\begin{eqnarray}\label{P_m expl.}
     &&P_m(x) = \sum\limits_{\ell=0}^{m-1}\sum\limits_{j=0}^{\ell}\sum\limits_{
          \alpha \in M^{\ell-j}}\frac{(-1)^{d\left(\alpha\right)}}{\alpha ! (j-\ell)!}
          \midint_{B_{\rho}}\left(D^\ell u\right)_{\rho}\left(y-x_0\right)^{\alpha}\, dy
          \left(x-x_0\right)^j,
\end{eqnarray}
where we use the following notation:
For $k \in \bN_0$ we denote by the number $d\left(\alpha \right) = k$ the order of the multi index $\alpha \in \bN^k$. In 
particular $d(0) = 0$. Furthermore $\left|\alpha\right| :=
\alpha_1 + \ldots + \alpha_k$ denotes the length of the multi index $\alpha$. 
The set
\begin{equation*}
     M^j = \left\{\alpha \in \left\{0\right\} \cup \bN^1 \cup \bN^2 \cup \ldots 
     \cup \bN^j:\ \left|\alpha\right| = j\right\}
\end{equation*}
denotes the set of all multi indices of order $\le j$ and length $j$. 
Moreover for a multi index $\alpha$ of order $k$ and $A \in \odot^\ell\left(\bR^n; 
\bR^N\right),\ell \ge k$ we set
\begin{eqnarray*}
     &&\midint_{\Omega}A\left(x-x_o\right)^{\alpha}\, dx = \midint_{\Omega}
          \ldots\midint_{\Omega}A\left(y_1-x_0\right)^{\alpha_1}\ldots
          \left(y_k-x_0\right)^{\alpha_k}\, dy_1\ldots dy_k.
\end{eqnarray*}
>From (\ref{P_m expl.}) we immediately get by differentiation for $k < m$:
\begin{equation*}
     D^kP_m(x) = \sum\limits_{\ell=k}^{m-1}\sum\limits_{j=k}^\ell\sum\limits_{
          \alpha\in M^{\ell-j}}\frac{(-1)^{d\left(\alpha\right)}j!}{\alpha ! (j-\ell)!
          (j-k)!}\midint_{B_{\rho}}\left(D^\ell u\right)_{\rho} \left(y-x_0 
          \right)^{\alpha}\, dy\left(x-x_0\right)^{j-k}.
\end{equation*}
Therefore we can estimate
\begin{eqnarray*}
     I_2
     & = &\midint_{B_{\rho}}\Biggl|\sum\limits_{\ell=k}^{m-1}\sum\limits_{j=k}^\ell\sum
          \limits_{\alpha \in M^{\ell-j}}\frac{(-1)^{d\left(\alpha\right)}j!}{
          \alpha ! (j-\ell)!(j-k)!}
          \midint_{B_{\rho}}\left(D^\ell u\right)_{\rho}\left(y-x_0 
          \right)^{\alpha}\, dy\left(x-x_0\right)^{j-k} - \left(D^k u\right)_{\rho} 
          \Biggr|^p\, dx\\
     & = &\midint_{B_{\rho}}\Biggl|\sum\limits_{\ell=k+1}^{m-1}\sum\limits_{j=k}^\ell\sum
          \limits_{\alpha \in M^{\ell-j}}\frac{(-1)^{d\left(\alpha\right)}j!}{
          \alpha ! (j-\ell)!(j-k)!}\midint_{B_{\rho}}\left(D^\ell u\right)_{\rho}\left( 
          y-x_0\right)^{\alpha}\, dy\left(x-x_0\right)^{j-k}\Biggr|^p\, dx\\
     &\le&c(p)\midint_{B_{\rho}}\sum\limits_{\ell=k+1}^{m-1}\sum\limits_{j=k}^\ell\sum
          \limits_{\alpha \in M^{\ell-j}}\frac{j!}{
          \alpha !(j-\ell)!(j-k)!}\left|\left(D^\ell u\right)_{\rho}\right|^p
          \midint_{B_{\rho}}
          \left|y-x_0\right|^{p\alpha}\, dy \left|x-x_0\right|^{p(j-k)}\, dx.\\
\end{eqnarray*}
where
\begin{eqnarray*}
     \midint_{B_{\rho}}\left|y-x_0\right|^{p\alpha}\, dx
     & = &\midint_{B_{\rho}}\ldots\midint_{B_{\rho}}\left|y_1-(x_0)_1\right|^{
          p\alpha_1}\ldots\left|y_{pd\left(\alpha
          \right)}-(x_0)_{d(\alpha)}\right|^{p\alpha_{d\left(\alpha\right)}}\, dy_1\ldots dy_{d
          \left(\alpha\right)}\\
     &\le&\rho^{p\alpha_1 + \ldots + p\alpha_{d\left(\alpha\right)}} 
          = \rho^{p\left|\alpha\right|} = \rho^{p(\ell-j)},
\end{eqnarray*}
when $\alpha \in M^{\ell-j}$.
Inserting this above we arrive at
\begin{equation*}
     I_2\le \sum\limits_{\ell=k+1}^{m-1}c(\ell,k)\rho^{p(\ell-k)}\left|\left(D^\ell u
          \right)_{\rho}\right|^p.
\end{equation*}
Combining the estimates for $I_1$ and $I_2$ we find that
\begin{equation*}
     \midint_{B_{\rho}} \left|D^k u(x) - \left(D^k u\right)_{\rho}\right|^p\, dx \le 
     c\left[\rho^{p(m-k)}\left(\left|D^m u\right|^p\right)_{\rho} + \sum\limits_{
     \ell=k+1}^{m-1} \rho^{p(\ell-k)}\left|\left(D^\ell u\right)_{\rho}\right|^p\right]
\end{equation*}
with $c \equiv c\left(n,p\right)$.
The claim follows now immediately by
\begin{eqnarray*}
     &&\midint_{B_{\rho}} \left|\delta u(x) - \left(\delta u\right)_{\rho}\right|^p
          \, dx \le c(n,m,p) \sum\limits_{k=0}^{m-1} \midint_{B_{\rho}}\left|D^k 
          u(x) - \left(D^k u\right)_{\rho}\right|^p\, dx.
\end{eqnarray*}
\end{proof}

\section{Proof of Theorem \ref{Hauptsatz}}

In the whole proof $u$ denotes the minimizer of the functional $\cF$ as required in 
Theorem \ref{Hauptsatz}. The proof will be divided in several lemmas. If the
function $u$ is mentioned there, we will always assume the minimizing property of $u$
without explicitely pointing it out again. 

\subsection{Localization, choice of radii and constants}\label{sec:funct.localization}

We start by Lemma \ref{High.Int.p(x)} which provides a higher integrability exponent $\delta_1$ such that for any $\Omega' \Subset \Omega$ there holds
\begin{equation*}
     \int_{\Omega'} \left|D^m u\right|^{p(x)\left( 1 + \delta_1\right)}
     \, dx < + \infty.
\end{equation*}
Since all our results are local we will assume for the rest of the proof that
\begin{equation*}
     \int_{\Omega} \left|D^m u\right|^{p(x)\left( 1 + \delta_1\right)}
     \, dx < + \infty.
\end{equation*}
Certainly we can choose $\delta_1$ so small such that
\begin{equation*}
     0 < \delta_1 \le \min \left\{\gamma_1 - 1,1\right\}.
\end{equation*}

Now let $1 < M < \infty$ and $\bar{L} \equiv \bar{L}(M)$ be the constant given by 
Lemma \ref{Techn.Lemma2}. We define the function 
\begin{equation*}
     g(z) := f_{U,A,\lambda} (z) = \lambda^{-2} \left[ f\left(x_0,U,A+\lambda z
     \right) - f\left(x_0,U,A\right) - \lambda Df\left(x_0,U,A\right)z\right],
\end{equation*}
with $U = (U_1,\ldots,U_m), |U_l| \le M \ (l=1,\ldots,m)$ and $|A| \le M$,
and find that by Lemma \ref{Techn.Lemma2} the function $g$ fulfills the
assumptions for Lemma \ref{High.Int.3} with $p \equiv p(x_0) = \const$. 
This yields a further higher integrability exponent $\delta_2 \equiv \delta_2(M)$. 
Subsequently we apply Lemma \ref{High.Int.Rand}, namely with the exponent
$q = p\left(1+\delta_1/4\right)$ and the constant $L$ replaced by $2L$ and therefore
obtain an up--to--the--boundary higher integrability exponent $\varepsilon$,
with $0 < \varepsilon < \delta_1/4$ and $\varepsilon \equiv \varepsilon(\gamma_1,
\gamma_2,L,\delta_1)$.

Now we set
\begin{equation}\label{Def.delta_3}
     \delta_3 := \min \left\{\varepsilon, \delta_2\right\} \equiv \delta_3(M)
\end{equation}
and choose a radius $R_M > 0$ such that
\begin{equation*}
     \omega_1\left(R_M\right) \le \frac{\delta_3}{4}.
\end{equation*}
From now on let $O \Subset \Omega$ be an open set whose diameter does not exceed  $R_M$.

On the set $O$ we can estimate the exponent $p$ at any point in terms of the maximal or
minimal exponent, respectively.
 
Therefore we set
\begin{equation}\label{Def.p2p1}
     p_1 := \inf\left\{p(x):\ x \in O \right\},\ \ \ p_2 := \sup\left\{
     p(x):\ x \in O\right\}.
\end{equation}
This implies (since $p_2 - p_1 \le \omega\left(\left|x - x'\right|\right) \le 
\omega\left(R_M\right) \le \delta_3 / 4 \le \varepsilon / 4 < \delta_1 / 16$):
\begin{eqnarray}
     p_2\left( 1 + \delta_1 / 4\right) &\le\ \ p_1 \left( 1 + \delta_1\right) 
          &\le\ p(x)\left( 1 + \delta_1\right),\label{HS:Absch.p1}\\
     p_2\left(1+\delta_2 / 4\right) &\le\ \ p_1\left( 1 + \delta_2\right) 
          &\le\ p(x) \left( 1 + \delta_2\right),\label{HS:Absch.p2}\\
     p_2\left( 1 + \varepsilon / 4\right) &\le\ \ p_1\left( 1 + \epsilon\right) 
          &\le\ p(x)\left( 1 + \varepsilon\right),\label{HS:Absch.p3}
\end{eqnarray}
for all $x \in O$.

\begin{remark}\label{Bem:epsilon}
$\delta_3$ does in fact depend on $M$ ( $\delta_3 \to 0$ as $M \to \infty$). 
$\varepsilon$ is independent of $M$ (see the dependencies of the exponent 
$\varepsilon$ in Lemma \ref{High.Int.Rand} for this purpose) and stays bounded away
from $0$ for any $M$.
\strut\hfill$\blacksquare$
\end{remark}
\begin{remark}\label{Bem:Lokalisierung}
By the constraint $\mbox{diam} (O) \le R_M$ the open set $O$ depends on the solution
$u$ of the variational problem itself. $O$ will be chosen in a special way at the end 
of the proof and it will be shown that the regular set is open.
\strut\hfill$\blacksquare$
\end{remark}

\subsection{Freezing}

Now we will
show that there exists a function $\breve{u}$ which is close to the
minimizer $u$ with respect to the $L^{p_2}$-- distance and which is an almost minimizer of the 'frozen' problem.

\begin{Lem}\label{freezing}
Let $M_1,M_2 > 0$ and $B\left(x_0,4R\right) \Subset O$. Then there exist $\beta_1,\beta_2$ dependent on $\gamma_1,\gamma_2,L,m$ and $\alpha$, but independent of
$M_1,M_2,R$ and $x_0$ and a constant $\check{c}(M_1,M_2)$ such that the following holds:
Let $u$ be a minimizer of the functional $\cF$ and $\left(|D^m u|^{p_2}\right)_{x_0,4R} 
\le M_1$, $\left|(D^k u)_{x_0,R}\right| \le M_2$, for $k=0,\ldots,m-1$, then there 
exists a function $\breve{u} \in u + W^{m,p_2}_0\left(B\left(
x_0,R\right);\bR^N\right)$ such that
\begin{eqnarray*}
     &(1)& \midint_{B\left(x_0,R\right)} \left|D^m u - D^m 
          \breve{u} \right|^{p_2} \, dx \le \check{c} R^{\beta_1},\\
     &(2)& \midint_{B\left(x_0,R\right)} f\left(x_0,\left(
          \delta u\right)_R,D^m \breve{u}\right)\, dx\\
     &&\mskip+30mu
          \le \midint_{B\left(x_0,R\right)}
          f\left(x_0,\left(\delta u\right)_R,D^m w\right)\, dx + R^{\beta_2} 
	  \midint_{B\left(x_0,R \right)} \left|D^m w - D^m 
	  \breve{u}\right|\, dx,\\
\end{eqnarray*}
for all $w \in u + W^{m,p_2}_0 \left(B\left(x_0,R\right); \bR^N\right).$ 
\strut\hfill$\blacksquare$
\end{Lem}

\begin{proof}
We consider the frozen integrand
\begin{equation*}
     g(z) := f\left(x_0,\left(u\right)_R,\left(Du\right)_R,\ldots,
          \left(D^{m-1}u\right)_R,z\right)
     \equiv f\left(x_0,\left(\delta u\right)_R,z\right)\\
\end{equation*}
and define $v \in u + W^{m,p\left(x_0\right)}_0\left(B_R;\bR^N\right)$ as
the unique solution of the minimization problem
\begin{equation*}
     \min\left\{ \int_{B_R} g\left(D^m w\right)\, dx :\ w \in u + 
     W^{m,p\left(x_0\right)}_0 \left(B_R; \bR^N\right)\right\}.
\end{equation*}
The existence of $v$ is guaranteed since the functional is quasiconvex.

\subsubsection*{HIGHER INTEGRABILITY}

Applying Lemma \ref{High.Int.p(x)} to the function $u$ leads to (using \eqref{HS:Absch.p1})
\begin{equation}\label{High.Int.u}
\midint_{B_{R}} \left|D^m u\right|^{p_2\left(1 + \delta_1/4
          \right)}\, dx \le c(M_1).
\end{equation}
Now, Lemma \ref{High.Int.Rand}, applied to $v$ with $g(x,z) 
\equiv g(z)$, $p \equiv p\left(x_0\right)$, $h \equiv u \in W^{m,q}\left(B_R
\right)$, $q \equiv p\left(1+\delta_1/4\right) > p$, $a(x) \equiv 1$ and 
$\tilde{m} = \delta_1/4$ provides $\varepsilon \in 
\left(0,\delta_1/4\right)$ such that
\begin{equation*}
\begin{aligned}
     \midint_{B_R} &\left|D^m v\right|^{p_2\left(1+\varepsilon/4\right)}
          \, dx \\
     &\le c\Biggl[ \left( \midint_{B_R} \left|D^m v\right|^{p(x_0)}\, dx
          \right)^{1+\epsilon}
          + \left(\ \midint_{B_R} \left|D^m u\right|^{
          p\left(x_0\right)\left(1+\delta_1/4
          \right)}\, dx\right)^{\frac{1+\varepsilon}{1+\delta_1/4}} + 1\Biggr]\\
     & = c \left[ (1) + (2) + 1 \right].
\end{aligned}
\end{equation*}
The estimate above gives for the second integral
\begin{equation*}
     (2) \le \left( \midint_{B_R} |D^m u|^{p_2}\, dx + 1 \right)^{\frac{1+\varepsilon}{1+\delta_1/4}}
     \le c(M_1).
\end{equation*}
For the first integral we use the minimizing property of $v$, combined with the
growth condition of $g$ obtaining
\begin{equation*}
\begin{aligned}
     \biggl(\midint_{B_R}&|D^m v|^{p(x_0)}\, dx\biggr)^{1+\varepsilon}\\
     & \le C \left(
     \midint_{B_R} |D^m u|^{p(x_0)}\, dx+1\right)^{1+\varepsilon} \le
     \left(\midint_{B_R} |D^m u|^{p_2}\, dx+1\right)^{1+\varepsilon} \le c(M_1).
\end{aligned}
\end{equation*}
Thus we have
\begin{equation}\label{High.Int.v}
\midint_{B_R}\left|D^m v\right|^{p_2\left(1+
     \epsilon/4\right)}\, dx \le c(M_1)
\end{equation}

\subsubsection*{A COMPARISON ESTIMATE}

We start by splitting as follows
\begin{eqnarray*}
     \midint_{B_R} \left[ g\left(D^m u\right) - g\left(D^m v\right)\right]
          \, dx
     &=& \midint_{B_R} \left[f\left(x_0,\left(\delta u\right)_R,D^m 
          u\right) - f\left(x_0,\left(\delta u\right)_R,D^m v\right)\right]
          \, dx\\
     &=& \midint_{B_R} \left[ f\left(x_0,\left(\delta u\right)_R,
          D^m u\right) - f\left(x_0,\delta u,D^m u\right)\right]
          \, dx\\
     &&\mskip+20mu + \midint_{B_R} \left[f\left(x_0,\delta u,D^m u\right) 
          - f\left(x, \delta u,D^m u\right)\right]\, dx\\
     &&\mskip+20mu + \midint_{B_R} \left[f\left(x,\delta u,D^m u\right) 
          - f\left(x,\delta u,D^m v\right)\right]\, dx\\
     &&\mskip+20mu + \midint_{B_R} \left[f\left(x,\delta u,D^m v\right) 
          - f\left(x_0,\delta u,D^m v\right)\right]\ dx\\
     &&\mskip+20mu + \midint_{B_R} \left[f\left(x_0,\delta u,D^m v\right) 
          - f\left(x_0,\left(\delta u\right)_R,D^m v\right)\right]\, dx\\
     &=& I^{(1)} + I^{(2)} + I^{(3)} + I^{(4)} + I^{(5)},
\end{eqnarray*}
with the obvious notation $I^{(1)}$ -- $I^{(5)}$. Subsequently we 
estimate $I^{(1)}$ -- $I^{(5)}$.

\noindent
{\bf Estimate for \boldmath $I^{(3)}$\unboldmath}: 
Using \eqref{Min.u} and \eqref{Stetigk.xi} we obtain:
\begin{eqnarray*}
     I^{(3)} 
     &\le& L\midint_{B_R}\omega_2\left(\left|\delta v - 
          \delta u\right|\right)\left(\mu^2 + \left|D^m v\right|^2
          \right)^{p_2/2}\, dx
          + L\midint_{B_R} \omega_2\left(\left|\delta v - 
          \delta u\right|\right)\, dx\\
     &=& I^{(3,1)} + I^{(3,2)}.
\end{eqnarray*}
We handle $I^{(3,1)}$ by the higher integrability result for $D^m v$:
Therefore let $r := p_2\left(1 + \tilde{\varepsilon}\right)$ with 
$\tilde{\varepsilon} \equiv \varepsilon/4$ \label{def.tildeeps}, and $\varepsilon$ the exponent of
(\ref{High.Int.v}). So $\tilde{\varepsilon} \in \left(0,\delta_1/4\right)$ and consequently $r \in (p_2, p_2(1+\delta_1/4))$.
Thus we obtain by H\"older's inequality, the boundedness of $\omega_2$, \eqref{High.Int.v},
Jensen's inequality and $\mu \le 1$:
\begin{eqnarray*}
     I^{(3,1)}
     &\le& c\left[\midint_{B_R} \left(\mu^2 + \left|D^m v(x)
          \right|^2\right)^{\frac{r}{2}}\, dx\right]^{\frac{p_2}{r}}
          \left[\midint_{B_R}\omega_2^{\frac{r}{r-p_2}}\left(
          \left|\delta v(x) - \delta u(x)\right|\right)\, dx\right]^{
          \frac{r-p_2}{r}}\\
     &\le&c\left[ \midint_{B_R} |D^m v|^{p_2\left(1+\varepsilon/4\right)}\, dx + 1\right]^{\frac{1}{
          1+\varepsilon/4}} \left[ \midint_{B_R} \omega_2\left( |\delta v - \delta u|\right)\, dx
          \right]^{\frac{r-p_2}{r}}\\
     &\le&c(M_1) \omega_2^{\sigma} \left( \midint_{B_R} |\delta v - \delta u|\, dx \right),
\end{eqnarray*}
where $\sigma \equiv \frac{r-p_2}{r}  = \frac{\tilde{
\varepsilon}}{1+\tilde{\varepsilon}}$. To estimate $I^{(3,2)}$ we use once again the boundedness and concavity of $\omega_2$ and Jensen's inequality, which together with the previous estimate for $I^{(3,1)}$ leads to
\begin{equation*}
     I^{(3)} \le c\omega_2^{\sigma}\left(\midint_{B_R} \left|\delta v - \delta u\right|\, dx \right),
\end{equation*}
with a constant $c \equiv c\left(\gamma_1,\gamma_2,L,n,m,M_1\right)$.

For estimating $\midint_{B_R}\left|\delta v - \delta u\right|\, dx$ we remark the
following: Since $u-v \in W^{m,p_2}_0\left( 
B_R\right)$, we can apply Poincar\'e's inequality obtaining
\begin{eqnarray*}
     \midint_{B_R} \left|D^{m-k}\left(u-v\right)\right|\, dx \le cR^k
     \midint_{B_R} \left|D^m\left(u-v\right)\right|\, dx,
\end{eqnarray*}
and therefore with \eqref{High.Int.u}, \eqref{High.Int.v} and $R \le 1$:
\begin{eqnarray*}
     \midint_{B_R}\left|\delta u - \delta v\right|\, dx 
     &\le&\sum_{k=0}^{m-1} \midint_{B_R}\left|D^k\left(u-v\right)\right|\, dx\\
     &\le&c\sum_{k=0}^{m-1} R^{m-k}\midint_{B_R} \left|D^m\left(u-
          v\right)\right|\, dx \le c\sum_{k=0}^{m-1} R^{m-k}\\
     &\le&c(n,m,M_1)R.
\end{eqnarray*}
Altogether we get
\begin{equation}
I^{(3)} \le c_{M_1}\omega_2^{\sigma}\left(R\right),
\end{equation}
with a constant $c_{M_1} \equiv c_{M_1}(n,m,L,M_1,\gamma_1,\gamma_2)$.

\noindent
{\bf Estimate for \boldmath $I^{(2)}$\unboldmath}:
We use \eqref{Stetigk.x} and \eqref{High.Int.u} to obtain
\begin{eqnarray*}
     &&I^{(2)} \le L\midint_{B_R} \omega_1\left(\left|x-x_0\right|\right)
          \left[\left(1+\left|D^m u\right|^2\right)^{p(x)/2} + \left(1+
          \left|D^m u\right|^2\right)^{p\left(x_0\right)/2}\right]\\
     &&\mskip+250mu \cdot
          \left[1+\log \left(1+\left|D^m u(x)\right|^2\right)\right]\, dx.
\end{eqnarray*}
We estimate the integrand in the following elementary way:
\begin{equation*}
     \left[\left(1+|z|^2\right)^{p\left(x_0\right)/2} + \left(1+|z|^2
          \right)^{p(x)/2}\right]\left(1+\log \left(1+|z|^2\right)\right)
     \le c(\gamma_1,\varepsilon) \left(1+|z|^{p_2(1+\varepsilon/4)}\right).
\end{equation*}
Thus we end up with
\begin{equation*}
     I^{(2)} \le c \omega_1(R) \left( \midint_{B_R} |D^m u|^{p_2(1+\varepsilon)}\, dx + 1 \right)
     \le c_{M_1} \omega_1(R).
\end{equation*}

\noindent
{\bf Estimate for \boldmath $I^{(1)}$\unboldmath}:
Using \eqref{Stetigk.xi}, the concavity of $\omega_2$ and \eqref{High.Int.u}, we 
obtain, proceeding in an analogue way to the estimate of $I^{(3)}$:
\begin{equation*}
     I^{(1)} 
      \le c_{M_1}\omega_2^{\sigma}\left(
          \midint_{B_R}\left|\delta u - \left(\delta u\right)_R\right|\, dx
          \right).
\end{equation*}
To estimate the integral we apply Lemma \ref{Niedr.Abl.}, obtaining by the
boundedness of the mean values
\begin{equation*}
     I^{(1)} 
      \le c_{M_1}\omega_2^{\sigma}\left(C\sum\limits_{k=0}^{m-1}
          \left[R^{m-k}c_{M_1} + \sum\limits_{l=k+1}^{m-1}cR^{l-k}c_{M_2}
          \right]\right),
\end{equation*}
which yields
\begin{eqnarray*}
     &&I^{(1)} \le c_{M_1,M_2}\omega_2^{\sigma}\left(R\right).
\end{eqnarray*}

\noindent
{\bf Estimate for \boldmath $I^{(4)}$\unboldmath}:
Completely analogous to the estimate of $I^{(2)}$ we use \eqref{Stetigk.x} and \eqref{High.Int.v} to find:
\begin{equation*}
     I^{(4)} 
     \le c\omega_1\left(R\right)\midint_{B_R}\left(1+\left|D^m v\right|^{p_2
          \left(1+\varepsilon/4\right)}\right)\, dx
     \le c_{M_1}\omega_1\left(R\right).
\end{equation*}

\noindent
{\bf Estimate for \boldmath$I^{(5)}$\unboldmath}:
This term is treated exactly as $I^{(1)}$ by
using \eqref{Stetigk.xi}, \eqref{High.Int.v} and the boundedness 
of the mean values $(D^k u)$:
\begin{equation*}
     I^{(5)}
     \le c_{M_1}\omega_2^{\sigma}\left(R\right).
\end{equation*}

Combining the estimates for $I^{(1)}$ -- $I^{(5)}$ we finally arrive at:
\begin{equation*}
     \midint_{B_R} \left[g\left(D^m u\right) - g\left(D^m v\right)\right]\, dx
          \le c_{M_1,M_2}\left(\omega_1\left(R\right) + \omega_2^{\sigma}\left(R\right)
          \right).
\end{equation*}
By assumption (\ref{Stet.Mod.stark}) we finally obtain for $R < 1$ and by letting $\tilde{\alpha} \equiv \sigma \alpha$:
\begin{equation}\label{Absch.u.v}
     \midint_{B_R} \left[g\left(D^m u\right) - g\left(D^m v\right)\right]\, dx 
     \le C\left(R^{\alpha} + R^{\sigma\alpha}\right) \le c_{M_1,M_2}R^{\tilde{
     \alpha}},
\end{equation}
with $\tilde{\alpha} \equiv \tilde{\alpha} \left( \alpha,\gamma_1,\gamma_2,L,m
\right)$ and the constant $c_{M_1,M_2}$ depending on $n,L,m,\gamma_1,\gamma_2, M_1, M_2$.

\begin{Bem}\label{Bem:alpha}
Since $\sigma \equiv \sigma\left(\varepsilon\right)$ is 
independent of $M_1,M_2$, also $\tilde{\alpha}$ does not depend on 
$M_1,M_2$ (see Remark \ref{Bem:epsilon} on page \pageref{Bem:epsilon} concerning the dependencies of $\varepsilon$).
\end{Bem}
\begin{Bem}
In estimate (\ref{Absch.u.v}) we used the boundedness of the mean values of $D^k u$
$(k=0,\ldots ,m-1)$ and of $|D^m u|^{p_2}$ on the balls $B_R$. However at the end of the 
proof we will define the regular set $\Omega_0$ in a way that these assumptions are 
satisfied automatically on $\Omega_0$.
\end{Bem}

\subsubsection*{A FURTHER COMPARISON FUNCTION}

Since the functional $\cF$ is only quasiconvex, we cannot directly estimate the $L^p$--distance of $D^m u$ and $D^m v$. Therefore we apply Ekeland's variational principle providing
a further function $\breve{u}$, which is close to the original minimizer $u$ with
respect to the $L^{p_2}$-- norm, therefore not anymore a minimizer of the frozen 
functional, but anyway an almost-- minimizer in the sense of \eqref{Fastmin.u}. 

We consider 
\begin{equation*}
\cX := u + W^{m,1}_0 \left(B_R; \bR^N\right),
\end{equation*}
together with
\begin{equation*}
     d: \cX \times \cX \to [0,\infty),\ \ \ \left(z_1,z_2\right) \mapsto 
     \hat{C}_M^{-1} R^{-\tilde{\alpha}/4}\midint_{B_R} \left|D^m z_1 - 
     D^m z_2 \right|\, dx.
\end{equation*}
Then the functional
\begin{equation*}
     {\mathcal G}: \cX \to \bR,\ \ \ {\mathcal G}(z) := \left\{ \begin{array}{ll} 
     \displaystyle{\midint_{B_R} g\left(D^m z\right)\, dx} & \text{ if } z 
     \in u + W^{m,p\left(x_0\right)}_0\left(B_R;\bR^N\right),\\ + \infty & 
     \text{ otherwise }.\end{array}\right.
\end{equation*} 
is obviously lower semicontinuous on the complete metric space $(\cX,d)$. By construction of $v$ and \eqref{Absch.u.v} we have
\begin{equation*}
     {\mathcal G}(v) = \min_X {\mathcal G}, \qquad {\mathcal G}(u) \le \inf_X {\mathcal G} + c_{M_1,M_2} R^{\tilde{\alpha}}.
\end{equation*}
Therefore Lemma \ref{eke_orig} provides $\breve{u} \in u + W^{m,p\left(x_0\right)}_0 \left(B_R;\bR^N\right)$ satisfying
\begin{eqnarray}
     \midint_{B_R} \left|D^m u - D^m \breve{u}\right|
          \, dx 
     &\le& c_{M_1,M_2} R^{\tilde{\alpha}/4} \ \ \text{ and }\label{Ek.1}\\
     \midint_{B_R}f\left(x_0,\left(\delta u\right)_R,D^m \breve{u}
          \right)\, dx
     &\le& \midint_{B_R} f\left(x_0,\left(\delta u
          \right)_R,D^m w\right)\, dx \nonumber\\
     &&\mskip+100mu + R^{\frac{3\tilde{\alpha}}{4}} \midint_{B_R} 
	  \left|D^m w - D^m \breve{u}\right|\, dx, \label{Fastmin.u}
\end{eqnarray}
for all $w \in u + W^{m,p(x_0)}_0 \left(B_R; \bR^N\right)$. This proves assertion $(2)$ with $\beta_2 := 3\tilde{\alpha}/4 \equiv \beta_2\left(L,m,\gamma_1,\gamma_2, 
\alpha\right)$.

For showing $(1)$ we consider the functional $\bar{\mathcal G}: u + W^{m,p(x_0)}_0(B_R; \bR^N) \to 
\bR$ defined by 
\begin{equation*}
     \bar{\mathcal G}(w) \equiv \midint_{B_R} f\left(x_0,\left(\delta u
          \right)_R, D^m w\right)\, dx + R^{3\tilde{\alpha}/4} \midint_{B_R} 
          \left|D^m w - D^m \breve{u} \right|\, dx.
\end{equation*}
By (\ref{Fastmin.u}) $\breve{u}$ minimizes $\bar{\mathcal G}$. Let
\begin{equation*}
     g(x,z) := f\left(x_0,\left(\delta u\right)_R,z\right) + R^{3\tilde{\alpha}/4} 
     \left|z - D^m \breve{u}(x)\right|.
\end{equation*}
It is easy to see that the growth condition \eqref{Wachst.2} translates into
\begin{equation}\label{growth.g}
     \tilde{L}^{-1} |z|^{p\left(x_0\right)} \le g\left(x,z\right) \le \tilde{L} 
     \left[ |z|^{p\left(x_0\right)} + a(x) \right],
\end{equation}
with $a(x) \equiv 1 + \left|D^m \breve{u}(x)\right| > 0$ and $\tilde{L} 
\equiv L + 1$. We apply Lemma \ref{High.Int.Rand} to the functional $\bar{\mathcal G}$ with $g(x,z) \equiv f
\left(x_0,\left(\delta u\right)_R,z\right) + R^{3\tilde{\alpha}/4}\left|z - 
D^m \breve{u}\right|$, $p\equiv p\left(x_0\right)$, $q \equiv p\left(x_0\right)
\left(1+\delta_1/4\right)$, $a \equiv \left|D^m \breve{u}\right| + 1 \in 
L^{\gamma_1}\left(B_R\right)$, $\gamma \equiv \gamma_1$, such that by $\tilde{m} 
\equiv \delta_1/4$ we obtain
\begin{eqnarray*}
     \midint_{B_R} \left|D^m \breve{u}\right|^{p_2\left(1+
          \varepsilon/4\right)}\, dx
     &\le&c\Biggl(\ \midint_{B_R} \left|D^m \breve{u}
          \right|^{p\left(x_0\right)}\, dx \Biggr)^{1+\varepsilon}
          + c\Biggl(\midint_{B_R}\left|D^m 
          u\right|^{p_2\left(1+\delta_1/4\right)}\, dx \Biggr)^{
          \frac{1+\varepsilon}{1+\delta_1/4}}\\
     &&\phantom{c\Biggl(\ \midint_{B_R} \left|D^m \breve{u}
          \right|^{p\left(x_0\right)}\, dx \Biggr)^{1+\varepsilon}} + c
          \Biggl(\ \midint_{B_R}
          \left|D^m \breve{u}\right|^{1+\delta_1/4}\, dx \Biggr)^{
          \frac{1+\varepsilon}{1+\delta_1/4}}\\
     &=&  c\left[ II^{(1)} + II^{(2)} + II^{(3)} \right].
\end{eqnarray*}
By \eqref{Fastmin.u}, \eqref{growth.g}, \eqref{Ek.1} and \eqref{High.Int.u} we obtain
\begin{equation*}
     II^{(1)} \le c_{M_1}. 
\end{equation*}
By higher integrability of $D^m u$, $II^{(2)}$ is also estimated by a constant $c_{M_1}$.
Finally H\"older's inequality leads to
\begin{equation*}
     II^{(3)} \le c \Biggl(\ \midint_{B_R} \left|D^m \breve{u}
     \right|^{p\left(x_0\right)}\, dx + 1\Biggr)^{1+\varepsilon} \le c_{M_1}.
\end{equation*}
Therefore we conclude
\begin{equation}\label{High.Int.breve u}
     \midint_{B_R} \left|D^m \breve{u}\right|^{p_2\left(1+\varepsilon/4
     \right)}\, dx \le c_{M_1}.
\end{equation}
Now we interpolate between $1$ and $p_2\left(1+\tilde{\varepsilon}\right)$ (with the definition
$\tilde{\varepsilon} = \varepsilon/4$ on page \pageref{def.tildeeps}).
For $ 0 < \theta < \frac{1}{p_2}$ we write by using H\"older's inequality 
\begin{eqnarray*}
     \midint_{B_R} \left|D^m u - D^m \breve{u}\right|^{p_2}\, dx
     &\le&\Biggl(\ \midint_{B_R} \left|D^m u - D^m \breve{u}
          \right|\, dx\Biggr)^{p_2 \theta} \Biggl(\ \midint_{B_R} \left|
          D^m u - D^m \breve{u}\right|^{\frac{p_2\left(1-\theta\right)}{1-p_2 
          \theta}}\, dx\Biggr)^{1-p_2\theta}\\
     &=&  II^{(4)} \cdot II^{(5)}.
\end{eqnarray*}
By (\ref{Ek.1}) we get for the first integral:
\begin{equation*}
     II^{(4)}
     \le \left(c_{M_1} R^{\tilde{\alpha}/4}\right)^{p_2\theta} = c_{M_1} R^{\theta
     \tilde{\alpha} p_2/4}.
\end{equation*}
Choosing $\theta=\frac{\tilde{ \varepsilon}}{p_2\left(1+\tilde{\varepsilon}\right)-1} 
< \frac{1}{p_2}$ and using higher integrability from \eqref{High.Int.breve u} and
\eqref{High.Int.u}
we have (note that $\tilde{\varepsilon} = \frac{\varepsilon}{4} 
\le \frac{\delta_1}{16} \le \frac{\delta_1}{4}$):
\begin{eqnarray*}
     II^{(5)} \le \Biggl(\midint_{B_R} \mskip-5mu \left|
          D^m u - D^m \breve{u}\right|^{p_2\left(1+\tilde{\varepsilon}\right)}\, dx
          \Biggr)^{\frac{1-\theta}{1+\tilde{\varepsilon}}} \le c_{M_1}.
\end{eqnarray*}
Thus we finally arrive at
\begin{equation*}
     \midint_{B_R} \left|D^m u - D^m \breve{u}\right|^{p_2}\, dx \le 
     \breve{C}_M R^{\theta \tilde{\alpha}/4}.
\end{equation*}
Noting that $\theta
     \ge \frac{\tilde{\epsilon}}{\gamma_2\left(1+\tilde{\epsilon}\right) - 1} 
     =: \bar{\theta}$
and $R \le 1$ we can estimate $R^{\tilde{\alpha}\theta/4}$ by $R^{\tilde{\alpha}\bar{\theta}/4}$.
Wo finish the proof we choose $\beta_1 \equiv \bar{\theta}\tilde{\alpha}/4$.
Note that by this choice $\beta_1 \equiv \beta_1(\gamma_1,\gamma_2,L,m,\alpha)$ is independent
of $M_1,M_2$, since $\varepsilon$ is so.
\end{proof}

\subsection{Excess--Decay estimate}

We set
\begin{equation}\label{Def.beta}
\beta := \frac{1}{2p_2}\min\left\{\beta_1,\beta_2\right\},
\end{equation}
and define the excess $\Phi(x_0,R)$ for all $x_0,R$ with $B\left(x_0,4R\right) \Subset O$ by letting
\begin{equation}\label{Def:Excess}
     \Phi\left(x_0,R\right) := \midint_{B\left(x_0,R\right)}
     \left|V_{p_2}\left(D^m u\right) - V_{p_2} \left(\left(D^m u\right)_{
     x_0,R}\right)\right|^2\, dx + R^{\beta},
\end{equation}
where $V_p$ is the function defined in (\ref{Def.V_p}) on page \pageref{Def.V_p} and $p_2$ is 
from \eqref{Def.p2p1}.

We define
\begin{equation}
     \qbar := \min \left\{2,p_2\right\},\ \ \ \ Q := \max\left\{2,p_2\right\}.
\end{equation}
\vspace{0.2cm}

\begin{Lem}[Excess decay]\label{Excess Decay}
Let $M>1$, $\beta$ from \eqref{Def.beta} and let $O \Subset \Omega$ be an open subset whose diameter does not exceed $R_M$, as explained on page \pageref{Def.delta_3}. Then there exist a constant 
$C_M$, depending also on $n,N,m,L,\gamma_1,\gamma_2$ and for every  
$\tau \in (0,1/24)$ a number 
$\varepsilon_0 \equiv \varepsilon_0\left(\tau,M\right)$, such that if
\begin{equation}\label{Vor.ED}
     \begin{array}{c}
     \left|\left(D^m u\right)_{x_0,\tau R}\right| \le M,\\[0.2cm]
     \left|\left(
          D^k u\right)_{x_0,R}\right| \le M\ \ \text{for}\ k=0,\ldots,m,\ \ \ \left|
          \left(D^m u\right)_{x_0,4R}\right| \le M, \\[0.2cm]
     \Phi\left(x_0,R\right) < \varepsilon_0,\ \ \ \ \Phi\left(x_0,4R\right) \le 1\\
     \end{array}
\end{equation}
hold on some $B(x_0,4R) \Subset \Omega$, then we have
\begin{equation}\label{Folge.ED}
     \Phi\left(x_0,\tau R\right) \le C_M \tau^{\beta} \Phi\left(x_0,R\right).
\end{equation}
\strut\hfill$\blacksquare$
\end{Lem}

\begin{proof}
We follow the ideas of \cite{Acerbi:2001}. Nevertheless
there are some modifications due to the higher order case which we will point out in the sequel.

For the whole proof we use the abbreviations $\Bh \equiv B(x_h,R_h)$ and $\Bvh \equiv
B(x_h,4R_h)$, as well as $(u)_{\Bh} \equiv (u)_{x_h,R_h}$, $(u)_{\Bvh} \equiv (u)_{x_h,
4R_h}$ respectively.

\noindent
{\bf STEP 1: BLOW UP:}
We prove the statement by contradiction. Therefore we 
assume that (\ref{Vor.ED}) holds, but (\ref{Folge.ED}) fails. Therefore there exists 
a sequence of balls $B\left(x_h,4R_h\right) \Subset O$, such that 
\begin{equation}\label{Vor.BU}
     \begin{array}{c}
     \left|\left(D^m u\right)_{x_h,\tau R_h}\right| \le M,\\
     \left|\left(
          D^k u\right)_{x_h,R_h}\right| \le M\ \text{for}\ k=0,\ldots,m,\ \ \
          \left|\left(D^m u\right)_{x_h,4R_h}\right| \le M,\\[0,2cm]
     \mu_h^2 := \Phi\left(x_h,R_h\right) \to 0\ \mbox{ as } h \to \infty,\ \ \ 
          \ \Phi\left(x_h,4R_h\right) \le 1,
     \end{array}
\end{equation}
but 
\begin{equation}\label{BU.schlecht}
     \Phi\left(x_h,\tau R_h\right) \ge C(M) \tau^{\beta} \Phi\left(x_h,R_h\right),
\end{equation}
where $C(M)$ will be chosen at the end of the proof. Without loss of generality
we can assume that $R_h \to 0$ as $h \to \infty$. Exactly as in the first order case
we see that there exists $c_M$ such that
\begin{equation}\label{BU.1}
     \left(\left|D^m u\right|^{p_2}\right)_{x_h,4R_h} \le c_M.
\end{equation}

By assumption we also have $|(D^k u)_{x_h,R_h}| \le M$ for $k=0,\ldots,m$.
We are now in a position to apply Lemma \ref{freezing} (with $M_1 := c_M$ from
\eqref{BU.1} and $M_2 = M$ from our hypothesis), which yields a sequence of
functions $u_h \in u + W^{m,p_2}_0\left(\Bh;\bR^N\right)$ 
satisfying
\begin{eqnarray}
     \mskip+30mu \midint_{\Bh}\left|D^m u - D^m u_h\right|^{p_2}\, dx 
     &\le& c_M R_h^{\beta_1},\label{u_h 1}\\
     \mskip+30mu \midint_{\Bh} f\left(x_h,U_h,D^m u_h\right) \, dx
     &\le& \midint_{\Bh} f\left(x_h,
          U_h,D^m w\right)\, dx + R_h^{\beta_2}
          \midint_{\Bh}\left|D^m w - D^m u_h
          \right|\, dx,\label{u_h 2}
\end{eqnarray}
for all $w \in u + W^{m,p_2}_0 \left(\Bh;\bR^N\right)$. We note that 
$\beta_1, \beta_2$ are independent of $h \in \bN$ and $M$. Here we used
the abbreviation $U_h \equiv \left(\delta u\right)_h \equiv \left(
\delta u\right)_{x_h,R_h}$. We let
\begin{equation*}
     A_h := \left(D^m u\right)_{\Bh},\qquad
     \lambda_h^2 := \midint_{\Bh}\left|
          V_{p_2}\left(D^m u_h\right) - V_{p_2}\left(A_h\right)\right|^2\, dx + 
          R_h^{\beta}.
\end{equation*}
We now rescale the functions $u_h$ in order to obtain a sequence of functions $v_h$ on the unit ball. Therefore let $P_h: \Bh \to 
\bR^N$ be the unique polynomial of degree $m$, for which there holds
\begin{equation}\label{Def.Ph}
     \int_{\Bh} D^k\left(u_h - P_h\right)\, dx 
     = 0\qquad (k = 0,\ldots,m-1), \qquad
     D^mP_h \equiv A_h.
\end{equation}
Now, using the notation $B_1 \equiv B(0,1)$ we define the sequence of rescaled functions by 
\begin{equation*}
     v_h\left(y\right) := \lambda_h^{-1}R_h^{-m}\left(u_h-P_h\right)\circ\left(
     x_h + R_hy\right),\ y \in B_1.
\end{equation*}
From \eqref{Def.Ph} we immediately see that
\begin{eqnarray}
     \left( D^k v_h\right)_{0,1} &=& 0 \quad (k=0,\ldots,m-1), \mbox{ and }
     \label{MW v_h}\\
     \quad D^m v_h\left(y\right) &=& \lambda_h^{-1}\left[D^m u_h \circ\left(x_h + R_hy
     \right) - A_h\right]. \label{Dm v_h}
\end{eqnarray}
Additionally we see, exactly as in the first order case, that
the sequences 
\begin{equation}\label{D^m v_h glm.beschr.}
\left(\left|D^m v_h\right|^{\qbar}\right)_{h \in \bN} \qquad \mbox{ and } \qquad \left(\lambda_h^{p_2-2} \left|D^m v_h\right|^{p_2}\right)_{h \in \bN} \mbox{ for } p_2 \ge 2
\end{equation}
are uniformly bounded in $L^1(B_1)$ by a constant $C_M$. This implies the existence 
of a subsequence --- without loss of generality the sequence $(v_h)$ itself --- and a
function $v \in W^{m,\qbar}\left(B_1;\bR^N\right)$ such that
\begin{itemize}\label{Konvergenz v_h}
\item[(a)] $v_h \to v$ weakly in $W^{m,\qbar}\left(B_1;\bR^N\right)$,
\item[(b)] $v_h \to v$ strongly in $W^{m-1,\qbar}\left(B_1;\bR^N\right)$,\\
     and therefore in particular $\left|v_h - v\right|^{\qbar} \to 0$ strongly in $L^1 
     \left(B_1\right)$,
\item[(c)] $\lambda_h^{p_2-2} \left|D^k \left(v_h-v\right)\right|^{p_2} \to 0$ 
     strongly in $L^1 \left(B_1\right), \ \ \text{for}\ \ k = 0,\ldots,m-1$, if $p_2 > 2$,
\item[(d)] $x_h \to x$ in $\bR^n$, with $x \in \bar{O}$,
\item[(e)] $A_h \to A$ in $\bR^{\cN}$, with $\left|A\right| \le M$.
\item[(f)] $U_h \to U$ in $\bR^{\cM}$, with 
     $U=(U_0,\ldots,U_{m-1})$ and $\left|U_k\right| \le M$ $(k=0,\ldots,m-1)$.
\end{itemize}
Using Lemma \ref{prop.V} and Jensen's inequality and proceeding exactly as in
\cite{Acerbi:2001} we obtain
\begin{equation}\label{lambda_h mu_h}
     \lambda_h^2 \le C_M \mu_h^2, \ \mbox{ and in particular }\ \lambda_h^2 \to 0.
\end{equation}

\noindent
{\bf STEP 2: $v$ SOLVES A LINEAR SYSTEM:}
The almost minimizing property of the functions $u_h$ directly translates into the following 
Euler--Lagrange system for the rescaled functions $v_h$:
\begin{equation}\label{Eulergl. 1}
     \midint_{B_1}\left< Df\left(x_h,U_h,A_h + 
     \lambda_h D^m v_h\right),D^m \phi\right>\, dy + \lambda_h I_h^{(2)} = 0,
\end{equation}
for all $\phi \in C_0^m\left(B_1,\bR^N\right)$, in which the second term satisfies
the estimate
\begin{equation}\label{Eulergl. 2}
     \lambda_h \left| I_h^{(2)}\right| \le R_h^{\beta_2} \midint_{B_1} 
     \left|D^m \phi\right|\, dy.
\end{equation}
Since $Df(x_h,U_h,A_h) = \const$, there holds additionally to the Euler--Lagrange
system:
\begin{eqnarray*}
     0 
     & = &\lambda_h^{-1} \int_{B_1} \bigl<Df\left(
          x_h,U_h,A_h+\lambda_h D^m v_h\right) - Df\left(x_h,U_h,A_h\right),D^m 
          \phi\bigr>\, dy + \left(II\right)_h\\
     &=: &I_h^{(1)} + I_h^{(2)},
\end{eqnarray*}
The definition of $\lambda_h^2$ and $\beta$ immediately imply
\begin{equation*}
     \lambda_h^{-1} R_h^{\beta_2} \to 0 \ \mbox{ and therefore }\ I_h^{(2)} \to 0.
\end{equation*}
For estimating $I_h^{(1)}$, we proceed exactly as in \cite{Acerbi:2001}, dividing the unit ball into
$E_h^+ := \{x \in B_1: \lambda_h\left|D^m v_h\right| \ge 1
\}$ and $E_h^- := B_1 \setminus E_h^+$, respectively the integral $I_h^{(1)}$ into 
$I_h^{(1),+}$ and $I_h^{(1),-}$. Deducing $ |E_h^+| \to 0$ by the uniform boundedness of 
$|D^m v_h|^{\qbar}$ in $L^1(B_1)$, and additionaly using the uniform boundedness
of $\lambda_h^{p_2-2}|D^m v_h|^{p_2}$ in $L^1(B_1)$ in the case $p_2 \ge 2$ together with the growth properties of $Df$ leads to 
\begin{equation*}
|I_h^{(1),+}| \to 0.
\end{equation*}
Using the fact that $f$ is of class $C^2$ with respect to $z$ and splitting the terms exactly as in \cite{Acerbi:2001}, we deduce (note that $|U_h| \le M$)
\begin{equation*}
     \lim_{h \to \infty} I_h^{(1),-} \to \int_{B_1} \left< D^2f(x,U,A)D^m v,D^m \varphi
     \right>\, dy.
\end{equation*}
Thus the function $v$ satisfies the following linear system:
\begin{equation}\label{v linear}
\int_{B_1} \left<D^2f\left(x,U,A\right)D^m v,D^m \phi\right>\, dy = 0,
\end{equation}
for all $\phi \in C_c^m\left(B_1;\bR^N\right)$. Furthermore the uniform strict
quasiconvexity of $f$ directly translates into the following property of $D^2 f$, which is equivalent to the Legendre--Hadamard condition:
\begin{equation}\label{Legendre-Hadamard}
     \int_O \left<D^2 f\left(x,U,A\right)D^m \phi,D^m \phi\right>\, dy \ge 
     c_M \int_O \left|D^m \phi\right|^2\, dy, 
\end{equation}
for all $\phi \in C_0^{\infty}\left(O;\bR^N\right)$ and with $0 \le c_M < 
\infty$.

Now the theory for linear elliptic systems of higher order with constant coefficients
(see \cite{Giaquinta:1983}) applies and yields $v \in C^{\infty}(B_1)$ together with the estimate
\begin{equation}\label{v glatt}
     \midint_{B_{\tau}} \left|D^m v - \left(D^m v\right)_{\tau}\right|^2\, dy
     \le c_M \tau^2,
\end{equation}
for all $\tau \le 1/4$.

\noindent
{\bf STEP 3: UPPER BOUND:}
We define the sequence of rescaled integrands
\begin{eqnarray*}
     f_h(z) &\equiv& f_{U_h,A_h,\lambda_h}(z)\\[0.1cm]
     &:= &\lambda_h^{-2} \left[ f\left(x_h,U_h,A_h+\lambda_h z
          \right) - f\left(x_h,U_h,A_h\right) - \lambda_h Df\left(x_h,U_h,A_h
          \right)z\right],
\end{eqnarray*}
for $h \in \bN,\ z \in \bR^{\cN}$.
For $r \in \left(0,1\right]$ and $w \in W^{m,1}\left(B_1;\bR^N\right)$ 
we define furthermore
\begin{equation*}
     I_h^r(w) := \int_{B_r} f_h\left(D^m w\right)\, dx.
\end{equation*}
We will show now, that for almost any $r \in (0, 1/6 )$ there holds
\begin{equation}\label{obere Schranke}
     \limsup\limits_{h \to \infty} \left[I_h^r\left(v_h\right) - I_h^r(v)
     \right] \le 0.
\end{equation}
Exactly as in the first order case we first observe that the minimality of $u_h$ 
translates into:
\begin{equation}\label{Min. I_h^r}
     I_h^r\left(v_h\right) \le I_h^r\left(v_h + \phi\right) + \lambda_h^{-2}
     R_h^{\beta_2}\int_{B_r}\left|D^m \phi\right|\, dy,
\end{equation}
for all $\phi \in W^{m,1}\left(B_1\right)$ with $\spt \phi \Subset B_r$.
Applying Lemma \ref{Techn.Lemma2} we observe that the hypothesis for
Lemma \ref{High.Int.3} are satisfied with $g \equiv f_h$, $\bar{u} \equiv v_h$ and
$\mu = \lambda_h^{-2}R_h^{\beta_2}$. With the choice of the quantities made at the
beginning of this section we observe that
\begin{equation}\label{High.Int.V D^m v_h}
     \midint_{B_{\frac{1}{6}}}\left|\frac{V_{p\left(x_h\right)}\left(\lambda_h
     D^m v_h\right)}{\lambda_h}\right|^{2\left(1+\delta_3\right)}\, dx \le c_M.
\end{equation}
In order to show \eqref{obere Schranke}, we proceed as in \cite{Acerbi:2001}: Consider the sequence of Radon measures $\alpha_h$ defined by
\begin{equation*}
     \alpha_h\left(A\right) := \int_A \lambda_h^{-2}\left[\left|V_{p_2}\left(
     \lambda_hD^m v_h\right)\right|^2 + \left|V_{p_2}\left(\lambda_h D^m v
     \right)\right|^2\right]\, dy.
\end{equation*}
for Borel sets $A \subset B_1$. Since $(\alpha_h)$ is uniformly bounded there exists a
subsequence -- again denoted by $(\alpha_h)$ -- and a Radon measure $\alpha$ with
$\alpha_h \overset{*}{\rightharpoonup} \alpha \ (h \to \infty)$. Since $\alpha(B_1) < \infty$, there
holds $\alpha(\partial B_t) = 0$ for all except at most countably many $t \in (0,1)$.
Thus without loss of generality we can assume $\alpha(\partial B_r) = 0$. Now let
$ 0 < s < r < 1$ and $\eta \in C_c^{\infty}\left(B_r\right)$ with
\begin{equation*}
     0 \le \eta \le 1,\ \ \eta \equiv 1 \ \text{ on }\ B_s,\ \ \left|D^k \eta
     \right| \le \frac{c(m)}{\left(r-s\right)^k}, \ \text{ for }\ k = 1,\ldots,m.
\end{equation*}
We test the minimality of $v_h$ in the sense of \eqref{Min. I_h^r} by the function
$\phi_h := \eta\left(v-v_h\right)$. Inserting the test function and using Lemma \ref{Techn.Lemma2}, a straight forward estimate under consideration of the properties of the function $V_p$ leads to
\begin{eqnarray*}
     I_h^r\left(v_h\right) - I_h^r\left(v\right)
     &\le&c \alpha_h\left(B_r \setminus B_s\right)\\
     &&\mskip+15mu
     + c\lambda_h^{-2}\sum\limits_{k=1}^m \frac{1}{
          \left(r-s\right)^{2kQ}}\int_{B_r \setminus B_s}\left|V_{p\left(
          x_h\right)}\left(\lambda_h D^{m-k}\left(v-v_h\right)\right)\right|^2
          \, dx + o_h\\
     &=&  (I)_h + (II)_h + o_h,
\end{eqnarray*} 
where $o_h$ denotes a quantity for which there holds $o_h \to 0$.
For proving that $(II)_h \to 0$, we use the strong convergence $\left|D^k\left(v-v_h\right)\right| \to 0$ in $L^{\qbar}(B_1)$
for all $k=0,\ldots,m-1$. In the case $p_2 \ge 2$, we do this by splitting
\begin{equation*}
     E^+_k :=\left\{x \in B_r \setminus B_s: \left|\lambda_h D^{m-k}\left(v-v_h
     \right)\right| \ge 1\right\},\quad E^-_k := (B_r \setminus B_s) \setminus E^+_k,
\end{equation*}
and using the properties of the function $V_p$ on the sets $E_k^+$ resp. $E_k^-$ to
estimate $(II)_h$ from above by terms containing $\| \lambda_h^{p_2-2}| 
D^{m-k} (v-v_h)|^{p_2}\|_{L^1(B_1)} $ and $\| |D^{m-k}(v-v_h)|^2\|_{L^1(B_1)}$. Due to the strong convergence we see that $(II)_h \to 0$ as $h \to \infty$. In the case
$1 < p_2 < 2$, we interpolate, defining $\theta := \frac{p(x_h)}{n+p(x_h)} < \frac{ 
1}{2},\  p^{\#} = p^{\#}(x_h) := \frac{2n}{n - p\left(x_h\right)}$ and estimating
by H\"older's inequality:
\begin{eqnarray*}
     &&\int_{B_r \setminus B_s}\left|V_p\left(\lambda_hD^{m-k}\left(
          v-v_h\right)\right)\right|^2\, dy\\
     &&\mskip+50mu \le \Biggl(\int_{B_r \setminus B_s}\mskip-5mu
          \left|V_p \left(\lambda_hD^{m-k}\left(v-v_h\right)\right)\right|\, dy
          \Biggr)^{2 \theta}\\
     &&\mskip+200mu \cdot
          \Biggl(\int_{B_r \setminus B_s}\mskip-5mu\left|V_p
          \left(\lambda_hD^{m-k}\left(v-v_h\right)\right)\right|^{p^{\#}}
          \, dy\Biggr)^{\frac{2\left(1-\theta\right)}{p^{\#}}}.
\end{eqnarray*}
While the first factor is simply estimated by $c(n,p)\lambda_h^{2\theta} ( \int_{B_r}
|D^{m-k}(v-v_h)|\, dx)^{2\theta}$, we use Sobolev--Poincar\'e's inequality in the
version of \cite{Duzaar:2004b}, Theorem 2 and Poincar\'e's inequality iteratively, to estimate the second factor from above, obtaining (note that $(\lambda_hD^{m-k}(v-v_h)
)_{0,1} = 0$)
\begin{eqnarray*}
     &&\Biggl(\int_{B_r \setminus B_s}\left|V_p
          \left(\lambda_hD^{m-k}\left(v-v_h\right)\right)\right|^{p^{\#}}
          \, dy\Biggr)^{\frac{2\left(1-\theta\right)}{p^{\#}}}\\
     &&\mskip+100mu \le c\Biggl(\int_{B_1}\left|V_p\left(\lambda_hD^{m-k+1}
          \left(v-v_h\right)\right)\right|^2\, dy\Biggr)^{1-\theta}\\
     &&\mskip+100mu \le c\biggl(\int_{B_1}\left|V_p\left(\lambda_hD^m 
          \left(v-v_h\right)\right)\right|^2\, dy \biggr)^{1-\theta}
          \le c_M \lambda_h^{2(1-\theta)}.
\end{eqnarray*} 
Taking these estimates together we conclude for any $p_2 > 1$:
\begin{equation*}
     I_h^r\left(v_h\right) - I_h^r\left(v\right)
     \le c_M \Biggl[ \alpha_h\left(B_r \setminus B_s\right) + 
          \sum\limits_{k
          =1}^m \frac{1}{\left(r-s\right)^{2kQ}}\left(
          \int_{B_r}\left|D^{m-k}\left(v-v_h\right)\right|\, dy\right)^{2\theta} 
          + o_h\Biggr],
\end{equation*}
with $c_M \equiv c_M(n,N,M,m,\gamma_1,\gamma_2,L,\nu)$. For $h \to \infty$ the right hand side converges to $\alpha\left(B_r \setminus 
B_s\right)$, since $\left|D^{m-k}\left(v-v_h\right)\right| \to 0$ strongly in $L^1$ for  
$k=1,\ldots,m$, and the assertion (\ref{obere Schranke}) follows with $s \to 
r$ since $\alpha\left(\partial B_r\right) =0$.

\noindent
{\bf STEP 4: LOWER BOUND:}
We will now show that 
\begin{equation}\label{untere Schranke}
     \limsup_{h \to \infty} \lambda_h^{-2} \int_{B_{r/2}}\left|V_{p_2}\left(
          \lambda_h\left(D^m v-D^m v_h\right)\right)\right|^2\, dy = 0,\ \ \text{for}\ 
          \ r \in \left(0,1/12\right).
\end{equation}
Therefore we proceed in a similar way to \cite{Acerbi:2001}: We consider $\frac{r}{2} < s < r < \frac{1}{12}$ and let
$\phi_h = \eta(v_h-v)$ with $\eta \in C^{\infty}_0(B_r), 0 \le \eta \le 1, \eta \equiv 1$ on $B_s$, $|D^k 
\eta| \le \frac{c}{(r-s)^k}$ for $k=1,\ldots,m$. Then we rewrite $I_h^r(v_h) - I_h^r(v)$ as follows:
\begin{equation*}
     I_h^r\left(v_h\right) - I_h^r\left(v\right) 
      = \left[I_h^r\left(v_h\right) - I_h^r\left(v+\phi_h\right)\right] + 
          \left[I_h^r\left(v+\phi_h\right) - I_h^r\left(v\right)\right]
      =:(I)_h + (II)_h.
\end{equation*}
In the sequel we estimate $(I)_h$ and $(II)_h$:

\noindent
{\bf Estimate for \boldmath $(I)_h$\unboldmath}: We proceed similarly to step 3 to 
conclude
\begin{equation*}
     \left|(I)_h\right| \le c_M\Biggl[ \alpha_h\left(B_r \setminus B_s\right) 
     + \sum\limits_{k=1}^m \frac{c}{\left(r-s\right)^{2kQ}}\Biggl(\ \int_{B_r}
     \left|D^{m-k}\left(v_h - v\right)\right|\, dy\Biggr)^{2\theta} + o_h
     \Biggr].\\
\end{equation*}

\noindent
{\bf Estimate for \boldmath $(II)_h$\unboldmath}: Here we follow the argumentation
of \cite{Acerbi:2001}. Therefore we only sketch the estimates here, pointing out
the additional difficulties due to the higher order case: First, we do the
splitting
\begin{eqnarray*}
     (II)_h 
     &=& \int_{B_r}\mskip-5mu f_h\left(D^m \phi_h\right)\, dy +
          \int_{B_r}\mskip-5mu \left[f_h\left(D^m v + D^m \phi_h\right) - 
          f_h\left(D^m v\right) - f_h\left(D^m \phi_h\right)\right]\, dy\\
     &=:&(II)_h^{(1)} + (II)_h^{(2)},
\end{eqnarray*}
with the obvious labelling.

\noindent
{\bf Estimate for \boldmath $(II)_h^{(1)}$\unboldmath}:
Lemma \ref{Techn.Lemma2} directly leads us to (note that $\eta \equiv 1$ on $B_{r/2}$)
\begin{equation*}
     (II)_h^{(1)} 
      \ge c_M^{-1}\lambda_h^{-2} \int_{B_{r/2}}\left|V_{p\left(x_h
          \right)}\left(\lambda_h\left(D^m v_h - D^m v\right)\right)\right|^2
          \, dy.
\end{equation*}

\noindent
{\bf Estimate for \boldmath $(II)_h^{(2)}$\unboldmath}: Exactly as in the first order
case we see that for any given $0 < \sigma \le 1$ there exists a ${\mathcal L}^n$--
measurable subset $S \subset B_1$ with 
\begin{equation*}
     {\mathcal L}^n\left(B_1 \setminus S\right) < \sigma \ \ \ \text{ and } \ \ \ 
     \lambda_h\left(\left|D^m v_h\right| + \left|D^m \phi_h\right|\right) \to 0 
     \ \text{ uniformly on } S.
\end{equation*}
Therefore we split as follows
\begin{eqnarray*}
     \left|(II)_h^{(2)}\right| 
     &\le&\int_{B_r \setminus S} \left| f_h\left(D^m\left(
          v + \phi_h\right)\right) - f_h\left(D^m v\right) - f_h\left(D^m \phi_h
          \right)\right|\, dy\\
     &&\mskip+40mu + \Biggl|\int_{B_r \cap S}\left(f_h\left(D^m\left(v + \phi_h
          \right)\right) - f_h\left(D^m v\right) - f_h\left(D^m \phi_h\right)
          \right)\, dy\Biggr|\\
     &=: &(II)_h^{(2,1)} + (II)_h^{(2,2)},
\end{eqnarray*}
with the obvious meaning of $(II)_h^{(2,1)}$ and $(II)_h^{(2,2)}$.

\noindent
{\bf Estimate for \boldmath $(II)_h^{(2,2)}$\unboldmath}: As in \cite[Proof of Prop. III.1, Step 5]{Acerbi:2001}, \cite[Proof of Prop. 3.4, Step 5]{Carozza:1998} we
see that
\begin{equation*}
     (II)_h^{(2,2)} \to 0 \ \ \ \mbox{ as }\ h \to \infty.
\end{equation*}

\noindent
{\bf Estimate for \boldmath $(II)_h^{(2,1)}$\unboldmath}: We set in the test function $\phi_h$ and proceed as in step 3 (using Lemma \ref{Techn.Lemma2}) to obtain
\begin{eqnarray*}
     (II)_h^{(2,1)} \mskip-10mu
     &\le&c\int_{B_r \setminus S}\left|
          \frac{V_{p\left(x_h\right)}\left(\lambda_hD^m v_h\right)}{\lambda_h}
          \right|^2\, dy + c\int_{B_r \setminus S}
          \left|\frac{V_{p\left(x_h\right)}\left(\lambda_h D^m v\right)}{
          \lambda_h}\right|^2\, dy\\
     &&\mskip+10mu + c\int_{B_r \setminus S}
          \lambda_h^{-2}\left|V_{p\left(x_h\right)}\left(\lambda_h \sum
          \limits_{k=1}^m \binom{m}{k} D^k\eta \odot D^{m-k}\left(v_h-v
          \right)\right)\right|^2\, dy.
\end{eqnarray*}
We use higher integrability \eqref{High.Int.V D^m v_h} together with H\"older's
inequality to estimate the first term above. The second term is controlled via the
smoothness ov $v$ on $B_1$, and the last term can be estimated exactly as in step 3 to conclude
\begin{equation*} 
     (II)_h^{(2,1)} \le c_M |B_r \setminus S|^{\frac{\delta_3}{1+
     \delta_3}} + \sum\limits_{k=1}^m \frac{C}{\left(r-s\right)^{2kQ}}\Biggl(
     \int_{B_r}\left|D^{m-k}\left(v-v_h\right)\right|\, dy\Biggr)^{2\theta} 
     + o_h.
\end{equation*}
Combining the above estimates, passing to the limit $h \to \infty$, taking into
account that $v_h \to v$ strongly in $W^{m-1,1}(B_1;\bR^N)$ and using the bound
>from above found in step 3, we conclude
\begin{equation*}
     \limsup_{h \to \infty} \lambda_h^{-2}\int_{B_{r/2}}\left|V_{p\left(x_h
          \right)}\left(\lambda_h\left(D^m v - D^m v_h\right)\right)\right|^2
          \, dy
     \le c_M\left[\alpha\left(B_r \setminus B_s\right) + \sigma^{
          \frac{\delta_3}{1+\delta_3}}\right].
\end{equation*}
Now we first let $\sigma \to 0$ and then $s \nearrow r$, obtaining
\begin{equation}\label{Sch. V_p x_h}
     \lambda_h^{-2}\int_{B_{r/2}} \left|V_{p\left(x_h\right)}\left(\lambda_h
     \left(D^m v - D^m v_h\right)\right)\right|^2\, dy = o_h,
\end{equation}
for almost every $0 < r < 1/12$. By the monotone dependence of the integral on the
domain of integration the statement holds for all $0 < r < 1/12$.

To finish the proof of step 4, we should carry over the 
above estimate for $V_{p(x_h)}$ to $V_{p_2}$. This is done in \cite[Proof of Prop. 4.2, step 5, p. 333]{Acerbi:2001} and applies for $D^m v$ instead of $Dv$ in
exactly the same way.

\noindent
{\bf STEP 5: FINISHING THE PROOF:}
Firstly, a straight forward calculation, using Lemma \ref{prop.V}, 
shows that by \eqref{u_h 1} there holds
\begin{equation}
     \mu_h^{-2} \midint_{B\left(x_h,\tau R_h\right)}\left|V_{p_2}\left(D^m u_h 
     - D^m u\right)\right|^2\, dx = o_h,
\end{equation}
and moreover
\begin{equation}
     \mu_h^{-2} \left|V_{p_2}\left(\left(
     D^m u_h\right)_{x_h,\tau R_h} - \left(D^m u\right)_{x_h,\tau R_h}\right)
     \right|^2 = o_h.
\end{equation}

We finish the proof by combining all the estimates we have shown before:
\begin{eqnarray*}
     &&\mskip-40mu\limsup_{h \to \infty} \mu_h^{-2} \Phi\left(x_h,\tau R_h\right) 
          \\
     &&\mskip+32mu \le \mskip+31mu c_M \limsup_{h \to \infty} \mu_h^{-2} 
          \midint_{ B\left(x_h,
          \tau R_h\right)} |V_{p_2}\left(D^m u - (D^m u)_{x_h,
          \tau R_h}\right)|^2\, dx\\
     &&\mskip+80mu + c_M \tau^{\beta} \limsup_{h \to \infty} \mu_h^{-2}
          R_h^{\beta}\\
     &&\mskip+19mu \overset{(\ref{lambda_h mu_h})}{\le} \mskip+18mu c_M \tau^{\beta}
          + c_M \limsup_{h \to \infty} 
          \mu_h^{-2} \midint_{B\left(x_h,\tau R_h\right)}|
          V_{p_2}\left(D^m u - D^m u_h\right)|^2\, dx\\
     &&\mskip+80mu + c_M \limsup_{h \to \infty} \mu_h^{-2} \midint_{B\left(x_h,
          \tau R_h\right)}|V_{p_2}\left(D^m u_h - (D^m u_h)_{ 
          x_h, \tau R_h}\right)|^2\, dx\\
     &&\mskip+80mu + c_M\limsup_{h \to \infty} \mu_h^{-2} |V_{p_2}\left( (D^m u_h
          )_{x_h,\tau R_h} - (D^m u)_{x_h,\tau R_h}\right)|^2\\
     &&\mskip+32mu = \mskip+31mu c_M \tau^{\beta}\\
     &&\mskip+80mu + c_M \limsup_{h \to \infty} 
          \mu_h^{-2} \midint_{B\left(x_h,\tau R_h\right)} |V_{p_2}
          \left(D^m u_h - (D^m u_h)_{x_h,\tau R_h}\right)|^2\, dx\\
     &&\mskip+0mu \overset{(\ref{Dm v_h}),(\ref{lambda_h mu_h})}{\le} c_M 
          \tau^{\beta} + c_M \limsup_{h \to \infty} \lambda_h^{-2} \midint_{B_{\tau}} \left|
          V_{p_2}\left(\lambda_h\left(D^m v_h - \left(D^m v_h\right)_{\tau}
          \right)\right)\right|^2\, dy\\
     &&\mskip+32mu \le \mskip+31mu c_M \tau^{\beta}
          + c_M \limsup_{h \to 
          \infty} \lambda_h^{-2} \midint_{B_{\tau}} \left|V_{p_2}\left(\lambda_h
          \left(D^m v_h - D^m v\right)\right)\right|^2\, dy\\
     &&\mskip+80mu + c_M \limsup_{h \to\infty} \lambda_h^{-2} \midint_{
          B_{\tau}} \left|V_{p_2}\left(\lambda_h\left(D^m v - \left(D^m v
          \right)_{\tau}\right)\right)\right|^2\, dy\\
     &&\mskip+80mu + c_M \limsup_{h \to \infty} \lambda_h^{-2}
          \left|V_{p_2}\left(\lambda_h\left(\left(D^m v
          \right)_{\tau} - \left(D^m v_h\right)_{\tau}\right)\right)\right|^2\\
     &&\mskip+0mu \overset{(\ref{v glatt}),(\ref{untere Schranke})}{\le} c_M
          \left(\tau^2 + \tau^{\beta}\right)\\[0.2cm]
     &&\mskip+32mu \le \mskip+31mu \hat{C}_M \tau^{\beta}.
\end{eqnarray*}
Choosing for example $C(M) := 2\hat{C}_M$, we 
obtain the desired contradiction. This completes the proof of the excess decay estimate.
\end{proof}

\subsection{Iteration of the excess decay estimate}

We will prove that under suitable smallness-- and boundedness conditions we can
iterate the excess decay estimate to obtain an estimate of the form
\begin{equation*}
     \Phi(x_0,\tau^k R) \le c \tau^{\alpha k}.
\end{equation*}
This is the tenor of the following
\begin{Lem}[Iteration]
Let $M \ge 2$, $B\left(x_0,16R\right) \Subset
 O \Subset \Omega$, where $O \equiv O_M$ is the open set from Lemma 
\ref{Excess Decay}. Let $c_M$ be the constant of Lemma \ref{Excess Decay} and let $0 < \tau < 1/24$ be chosen in such a way that 
$c_M \tau^{\beta/2} < 1/4$. Then there exists $\eta \equiv \eta\left(M,
\tau\right) \equiv \eta\left(M\right) \le \varepsilon_0\le 1$, where
$\varepsilon_0$ is the constant from Lemma \ref{Excess Decay}, such that if
\begin{eqnarray}\label{Vor.Iteration.1}
     &\left|\left(D^m u\right)_{x_0,\tau R}\right| + \left|\left(D^m u
          \right)_{x_0,R}\right| + \left|\left(D^m u\right)_{x_0,4R}\right| 
          \le M/4& \text{ and }\\[0.1cm]
     &\left|\left(D^j u\right)_{x_0,R}\right| \le M/4 \ \ \ \text{for}
          \ \ j=0,\ldots,m-1,\ , \nonumber
\end{eqnarray}
and
\begin{equation}\label{Vor.Iteration.2}
     \Phi\left(x_0,R\right) \le \eta,\ \ \Phi\left(x_0,4R\right) \le 1
\end{equation}
are satisfied, then there holds 
\begin{equation}
     \left|\left(D^m u\right)_{x_0,\tau^kR}\right| \le M,\ \ \ \Phi\left(x_0,
     \tau^kR\right) \le \tau^{k\beta/2} \ \ \ \text{for all}\ \ k \ge 1.
\end{equation}
\strut\hfill$\blacksquare$
\end{Lem}

\begin{proof}
The proof of this lemma is in many points identical with the proof in the first order case (see for example \cite{Fusco:1985} for the iteration scheme). The only additional
difficulty consists of veryfying $(\ref{Vor.Iteration.1})_2$ in every step. Therefore
we only point out how to do this. 
So we assume that for $k = 0,\ldots,s$ there holds
\begin{equation*}
     \Phi(x_0,\tau^k R) \le (C_M \tau^{\beta/2})^k \Phi(x_0,R), \qquad
     \Phi(x_0,4\tau^kR) \le c(M) \tau^{-2n}\Phi(x_0,\tau^{k-1}R).
\end{equation*}
To show $(\ref{Vor.Iteration.1})_2$ in the case $p_2 \ge 2$ we start by writing
\begin{eqnarray*}
     \left|(D^j u)_{\tau^s R} \right| 
     &\le& \left|(D^j u)_R \right| + \sum\limits_{k=0}^{s-1}\left|(D^j u)_{\tau^k R} -
           (D^j u)_{\tau^{k+1}R} \right|\\
     &\le& \frac{M}{4} + \sum\limits_{k=0}^{s-1} \left( \midint_{B_{\tau^{k+1}R}}
           \left| D^j u - (D^j u)_{\tau^k R} \right|^2\, dx \right)^{1/2}\\
     & = & \frac{M}{4} + \sum\limits_{k=0}^{s-1} (II)_k,
\end{eqnarray*}
with the obvious meaning of $(II)_k$. To estimate $(II)_k$, we note that $(DD^j u)_{
\tau^kR}(x-x_0)$ has meanvalue zero on balls centered in $x_0$. Therefore we obtain, applying
Poincar\'e's inequality iteratively:
\begin{eqnarray*}
     (II)_k 
     &\le& \tau^{-n/2} \left[ \midint_{B_{\tau^k R}} \left| D^j u - (D^j u)_{\tau^k R} 
           \right|^2\, dx \right]^{1/2}\\
     & = & \tau^{-n/2} \left[ \midint_{B_{\tau^k R}} \left| D^j u - (D^j u)_{\tau^k R}
           - (D D^j u)_{\tau^k R}(x-x_0)\right|^2\, dx \right]^{1/2}\\
     &\le& \tau^{-n/2} \left[ c_P \tau^k R \left( \midint_{B_{\tau^k R}} \left|DD^j u -
           (D D^j u)_{\tau^k R}\right|^2\, dx\right)^{1/2} \right]\\
     &\le& \tau^{-n/2} \Biggl[ c_P^2 (\tau^k R)^2 \left( \midint_{B_{\tau^k R}} \left|
           D^{j+2} u - (D^{j+2}u)_{\tau^k R}\right|^2\, dx \right)^{1/2}\\
     &&\mskip+280mu + c_P \tau^k R \left| (D^{j+2} u)_{\tau^k R}\right|\Biggr] \le \ldots \\
     &\le& \tau^{-n/2} \Biggl[ \left(c_P\tau^k R\right)^{m-j} \left( \midint_{B_{
           \tau^k R}} \left|D^m u - (D^m u)_{\tau^k R}\right|^2\, dx \right)^{1/2}\\
     &&\mskip+280mu +  \sum\limits_{l=2}^{m-j} c_P^{l-1} \left(\tau^k R
           \right)^{l-1} \left|(D^{j+l}u)_{\tau^k R}\right| \Biggr].
\end{eqnarray*}
Since by the induction hypothesis $(D^l u)_{\tau^k R} \le M$ for $k=0,\ldots,s-1$ and $l=0,\ldots,m$, by the definition
of the excess and with $\tau^k R \le 1$ and Lemma \ref{prop.V} we obtain
\begin{eqnarray*}
     (II)_k 
     &\le& \tau^{-n/2} \left[ (c_P \tau^k R)^{m-j} \Phi^{1/2}(x_0,\tau^k R) +
           \Phi(x_0,\tau^k R) M\sum\limits_{l=2}^{m-j} c_P^{l-1} (\tau^k R)^{l-2}
           \right]\\
     &\le& \tau^{-n/2} c_P^{m-j} \Phi^{1/2}(x_0,\tau^k R) + c_M\tau^{-n/2}\Phi(x_0,\tau^k R).
\end{eqnarray*}
Thus we end up with
\begin{equation*}
     \left| (D^j u)_{\tau^{s}R} \right| \le \frac{M}{4} + \tau^{-n/2} 
     \sum\limits_{k=0}^{s-1}\left[ c_P^{m-j} \Phi^{1/2}(x_0,\tau^k R) + c\Phi(x_0,\tau^k R) 
     \right].
\end{equation*}
The induction hypothesis together with $c_M\tau^{\beta/2} \le 1$ leads to
\begin{eqnarray*}
     \left| (D^j u)_{\tau^{s}R}\right| 
     &\le& \frac{M}{4} + \tau^{-n/2} \sum\limits_{k=0}^{s-1} \left[ c_P^{m-j} 
           \left(c_M
           \tau^{\beta/2}\right)^{k/2}\Phi^{1/2}(x_0,R) + c\left(c_M \tau^{\beta/2}
           \right)^k \Phi(x_0,R) \right]\\
     &\le& \frac{M}{4} + c\tau^{-n/2} \frac{\Phi^{1/2}(x_0,R)}{1- \sqrt{c_M\tau^{
           \beta/2}}} \le M,
\end{eqnarray*}
provided that $\eta \le \frac{9M^2\tau^n}{16c^2}\left(1-\sqrt{c_M\tau^{\beta/2}}\right)^2$. 
Note that in the case $1 < p_2 < 2$ we can estimate similary
to the case $p_2 \ge 2$, using the properties of the function $V_p$, and ending up with
the same smallness condition on $\eta$ (eventually with a modified constant $c_M$).
\end{proof}

\subsection{Construction of the regular set}

We will show that the hypothesis of the excess decay
estimate are satisfied on an open set $\Omega_0$ of full $n$ dimensional
Lebesgue measure.

First we define
\begin{equation}\label{Reg.M.}
\begin{aligned}
     \Sigma_1 &:= \left\{x_0 \in \Omega:\ \limsup\limits_{\rho \searrow 0} \left(
          \left|D^m u\right|^{p(x)}\right)_{x_0,\rho} = + \infty\right\},\\
     \Sigma_2 &:= \left\{x_0 \in \Omega:\ \limsup\limits_{\rho \searrow 0} \left|
          \left(D^k u\right)_{x_0,\rho}\right| = + \infty \ \text{ for }\ k=0,
          \ldots,m-1\right\},\\
     \Sigma_3 &:= \left\{x_0 \in \Omega:\ \limsup\limits_{\rho \searrow 0} \midint_{
          B\left(x_0,\rho\right)}\left|D^m u - \left(D^m u\right)_{x_0,\rho}
          \right|\, dx > 0\right\}.
\end{aligned}
\end{equation}
We will show that the smallness conditions from Lemma \ref{Excess Decay} are fulfilled on the set
\begin{equation*}
     \Omega_0 \equiv \Omega \setminus \left( \Sigma_1 \cup \Sigma_2 \cup \Sigma_3\right).
\end{equation*}
\begin{Bem}
The set $\Omega_0$ is of full Lebesgue measure, which can be directly seen by 
Lebesgue's theorem, since by $u \in W^{m,1}_{loc}(\Omega)$ there holds
$D^k u \in L^1_{loc}(\Omega)$ for $k=0,\ldots,m$.
\end{Bem}
Let $x_0 \in \Omega_0$, i.e.
\begin{eqnarray}
     &&\limsup_{\rho \searrow 0} \left( |D^m u|^{p(\cdot)} \right)_{x_0,\rho} < + \infty, \label{REG 1}\\
     &&\limsup_{\rho \searrow 0} \left| (D^k u)_{x_0,\rho}\right| < + \infty,\quad \text{for} \quad k=0,\ldots
       ,m-1 \label{REG 2}\\
     &&\liminf_{\rho \searrow 0} \midint_{B(x_0,\rho)} \left| D^m u - (D^m u)_{x_0,\rho}\right|\, dx = 0. 
       \label{REG 3}
\end{eqnarray}
Let $c_0 > 1$ be the constant of the higher integrability Lemma \ref{High.Int.p(x)}.
With \eqref{REG 1} and \eqref{REG 2} we can assume that there exists $M \ge \max\{2,8c_0\}$ such that
\begin{equation}\label{Fix.M}
     \left(\left|D^m u\right|^{p(x)}\right)_{x_0,2\rho} < \left( \frac{M}{8c_0}
          \right)^{\frac{1}{1+\delta_1}},\ \ \left|\left(D^k u\right)_{x_0,\rho}
          \right| \le \frac{M}{4} \quad \text{for} \quad k=0,\ldots,m\ ,
\end{equation}
for all radii $0 < \rho < \rho_0$ with $\rho_0 > 0$. With \eqref{REG 3} there exists a sequence of radii
$\rho_h \downarrow 0$ such that
\begin{equation*}
     \midint_{B(x_0,\rho_h)} \left|D^m u - (D^m u)_{x_0,\rho_h}\right|\, dx \to 0 \quad \text{as} \quad
     h \to \infty.
\end{equation*}
We consider such a sequence $\rho_h$. Then a straight forward estimate shows
\begin{equation*}
     \midint_{B(x_0,\rho_h/4)} \left| D^m u - (D^m u)_{x_0,\rho_h/4} \right| \, dx
     \le 2\cdot 4^n \midint_{B(x_0,\rho_h)} \left|D^m u - (D^m u)_{x_0,\rho_h} \right|\, dx,
\end{equation*}
and
\begin{equation*}
     \midint_{B(x_0, \tau \rho_h/4)} \left|D^m u - (D^m u)_{x_0,\tau \rho_h/4} \right|\, dx
     \le 2\left(4/\tau\right)^n \midint_{B(x_0,\rho_h)} \left|D^m u - (D^m u)_{x_0,\rho_h}\right|\, dx.
\end{equation*}
The constant $M$ from above fixes $c_M$ of the excess decay Lemma \ref{Excess Decay}. Therefore $\tau \equiv \tau(M)$ in
the excess decay estimate is fixed, and this again fixes the smallness parameter $\eta \equiv 
\eta(M,\tau)$. 

Let $B<1$ and $\theta$ be chosen in such a way that there holds
\begin{equation*}
     2^{\left(\gamma_2\right)^2}M\tilde{c}_M\left(B^{p_2\theta} + B^{2\theta}
     \right) \le \frac{\eta}{4} \ \ \ \text{ and }\ \ \ \ 1 = p_2\theta + 
     \frac{1-\theta}{1 + \delta_1/4},
\end{equation*}
where $\tilde{c}_M$ denotes the square of the constant appearing in Lemma \ref{prop.V}. So we have
\begin{equation*}
     \frac{1}{p_2} > \theta = \frac{\delta_1/4}{p_2\left(1+\delta_1/4\right) 
     -1} \equiv \theta\left(\delta_1,p_2\right)> 0,\ \ \ \ B \equiv B\left(M,
     \tilde{c}_M,p_2,\theta,\eta\right).
\end{equation*}
Then with the considerations above with $\rho_h =: 4R$ (for some $h \gg 1$) we have
\begin{equation}\label{Excess.smallnessB}
     \midint_{B\left(x_0,\sigma\right)}|D^m u - (D^m u)_{x_0,\rho}|\, dx < B,
\end{equation}
for $\sigma = 4R, \sigma = R$ and $\sigma = \tau R$.
Moreover we assume that $R$ is so small that
\begin{equation*}
     \left(4R\right)^{\beta} < \frac{\eta}{4}.
\end{equation*}
Then by \eqref{HS:Absch.p1} and the
higher integrability Lemma \ref{High.Int.p(x)} one can easily see
\begin{equation*}
     \left(\left|D^m u\right|^{p_2\left(1 + \delta_1/4\right)}\right)_{x_0,
          \sigma} 
     \le \frac{M}{2}.
\end{equation*}
Furthermore by interpolation between $1$ and $p_2\left(1+\delta_1/4\right)$ one can see
\begin{equation*}
     \midint_{B\left(x_0,\sigma\right)} |D^m u - (D^m u)_{x_0,
          \sigma}|^{p_2}\, dx
     \le 2^{\left(\gamma_2\right)^2}MB^{\theta p_2}.
\end{equation*}
Consequently we can estimate the excess $\Phi(x_0,\sigma)$ as follows: In the case
$1 < p_2 < 2$ by $|(D^m u)_{x_0,\sigma}| \le M$ we obtain with
the properties of the function $V_p$:
\begin{equation*}
     \Phi\left(x_0,\sigma\right)
      \le c_M \midint_{B\left(x_0,\sigma\right)}
          |D^m u - (D^m u)_{x_0,\sigma}|^{p_2}\, dx + 
          \sigma^{\beta}\\
     \le 2^{\left(\gamma_2\right)^2}MB^{\theta p_2} + 
          \sigma^{\beta} \le \eta.
\end{equation*}
In the case $p_2 \ge 2$ we obtain by H\"older's inequality and the properties of
$V_p$:
\begin{eqnarray*}
     \Phi\left(x_0,\sigma\right) \mskip-5mu
     &\le&\tilde{c}_M \Biggl[\midint_{B\left(
          x_0,\sigma\right)}\mskip-5mu |D^m u - (D^m u)_{x_0,\sigma}
          |^{p_2}\, dx
          + \Biggl(\midint_{B\left(x_0,\sigma
          \right)}\mskip-5mu |D^m u - (D^m u)_{x_0,\sigma}|^{p_2}
          \, dx\Biggr)^{\frac{2}{p_2}}\Biggr] + \sigma^{\beta}\\
     &\le&\hat{c}_M\Bigl[ 2^{\left(\gamma_2\right)^2}MB^{\theta p_2} 
          + 2^{\left(2/p_2\gamma_2\right)^2}M^{2/p_2}B^{2\theta}\Bigr] + 
          \sigma^{\beta}\\
     &\le&\hat{c}_M M 2^{\left(\gamma_2\right)^2}\left[B^{\theta p_2} 
          + B^{2\theta}\right] + \left(4R\right)^{\beta} \le \frac{\eta}{2} 
          \le \eta.
\end{eqnarray*} 
Thus the conditions for the excess decay estimate are satisfied.

\subsection{Localization and Conclusion}\label{sec:Lokalisierung}

Now we will show that 
$\Omega_0$ is actually open and that for $x_0 \in \Omega_0$, if 
(\ref{Vor.Iteration.1}) and (\ref{Vor.Iteration.2}) hold for a suitable $M$, 
there holds
\begin{equation}\label{Campanato}
     \midint_{B\left(x_0,\rho\right)}\left|V_{p_2}\left(D^m u\right) - \left(
     V_{p_2}\left(D^m u\right)\right)_{x_0,\rho}\right|\, dx \le C_M \rho^{
     \beta/4},
\end{equation}
for all $0 < \rho \le \tilde{R}_M$.

Therefore, let $x_0 \in \Omega_0$ and $M$ be fixed as in \eqref{Fix.M}.
Moreover let $R_M > 0$ be chosen as in the beginning, i.e.
\begin{equation*}
     \omega\left(R_M\right) \le \frac{\delta_3}{4} \ \ \ \text{ with }
     \ \ \ \delta_3 \equiv \delta_3(M),
\end{equation*}
where $\delta_3$ denotes the higher integrability exponent from \eqref{Def.delta_3}.
Let furthermore be $R < R_M/32$, such that (\ref{Excess.smallnessB}) holds for $\sigma = 
\tau R$, $\sigma = R$ and $\sigma = 4R$.

We set
\begin{equation}\label{Wahl O_M}
     O_M := B\left(x_0,R_M\right).
\end{equation}
By step 2 the conditions (\ref{Vor.Iteration.1}) and (\ref{Vor.Iteration.2}) for 
the iteration are fulfilled. The iteration (step 1) provides
\begin{equation}\label{nach Iteration}
     \left|\left(D^m u\right)_{x_0,\tau^k R}\right| \le M,\ \ \ \Phi\left(x_0,
     \tau^k R\right) \le \tau^{k \beta/2},
\end{equation}
for all $k \ge 1$.
Thus statement (\ref{Campanato}) follows immediately in the case $\rho = \tau^k R$ 
since
\begin{equation*}
     \midint_{B\left(x_0,\rho\right)}\mskip-20mu |V_{p_2}(D^m u) - \left(
          V_{p_2}(D^m u)\right)_{x_0,\rho}|^2\, dx
     \le \Phi\left(x_0,\rho\right)
     \le C_M \rho^{\beta/4}.
\end{equation*}
The statement for arbitrary radii we obtain by interpolation.

Now by construction of $O_M$ for any point $x_1 \in \Omega_0 \cap
B\left(x_0,R\right)$ with $R \le R_M/32$ there holds
\begin{equation}
     B\left(x_1,R\right) \subset B\left(x_0,R_M\right) = O_M \ \ \ \
     \mbox{ for all }\ R \le R_M/32.
\end{equation}
Thus also on $B\left(x_1,R\right)$ for $R \le R_M/32$ the conditions for the
iteration of the excess function are satisfied. In conclusion the regular set 
$\Omega_0$ is open.

Now we have shown that
\begin{equation}
     \midint_{B(x,\rho)} |V_{p_2}(D^m u) - \left( V_{p_2} (D^m u) 
     \right)_{x,\rho}|\, dy \le C_M \rho^{\beta/4},
\end{equation}
for all $x \in B\left(x_0,R_M/32\right)$ and for any radius $0 < \rho
< R_M/32$. Campanato's integral characterization implies that $V_{p_2}(D^m u)$ is H\"older continuous with 
exponent $\beta/4$ on $B\left(x_0,R_M/32\right)$.

Now, as one can easily see, H\"older continuity carries over from $V_{p_2}(D^m u)$ to $D^m u$ itself,
with exponent $\tilde{\beta} \equiv \min\{1,2/p_2\} \beta/8$.

By (\ref{Def.beta}) and (\ref{Def.p2p1}) the H\"older exponent $\tilde{\beta}$ 
depends on the local situation (in particular on $O_M$). However the exponents
$\beta_1$ and $\beta_2$ only depend on the global bounds $\gamma_1,\gamma_2$ of the
function $p$ (see Lemma \ref{freezing}). More precisely by
\begin{equation*}
     \tilde{\beta} = \min\{1,2/p_2\} \beta/8,\ \ \ \ \ \beta=1/2p_2 
     \min\{\beta_1,\beta_2\}
\end{equation*}
there holds
\begin{equation}
     \tilde{\beta} \in \left[ \tilde{\beta_1},\ \tilde{\beta_2} \right],
\end{equation}
in which
\begin{eqnarray*}
     \tilde{\beta_1} 
     &\equiv& \min\{1,2/\gamma_2\} 1/(16\gamma_2) \min\{\beta_1,
     \beta_2\},\\ 
     \tilde{\beta_2}
     &\equiv& \min\{1,2/\gamma_1\} 1/(16\gamma_1)\min\{\beta_1, \beta_2\}.
\end{eqnarray*}
By a covering argument we can show H\"older continuity of $D^m u$ on the whole
set $\Omega_0$.

Thus the proof is finished. \hfill \boxed{QED}


\bibliographystyle{plain}
\label{bibliography}
\makeatletter
\addcontentsline{toc}{chapter}{\bibname}
\makeatother
\bibliography{references}
\nocite{Acerbi:2001b,Acerbi:2002,Coscia:1999,Eleuteri:2004,Evans:1986,Evans:1987,Kristensen:2006,Kronz:2002,Manfredi:1988,Marcellini:1989,Rajagopal:2001,Ruzicka:2000}

\end{document}